\documentclass[global,final]{svjour}[21.1.2005]
\journalname{Combinatorica}

\usepackage {latexsym} \usepackage {amsmath} \usepackage {amsfonts}
\usepackage {amssymb}
\usepackage {times}
\usepackage {epsf}
\usepackage [PostScript=dvips]{diagrams}
\usepackage {paralist}
\smartqed
\spnewtheorem*{problem*}{Problem}{\bfseries}{\rmfamily}

\let\phi\varphi
\def\?#1{\par{\bf ??? #1 !!!}\par}
\def\sss{\scriptscriptstyle}
\def\zet{{\mathbb Z}}
\def\en{{\mathbb N}}
\def\Cay{\mathop{\rm Cay}\nolimits}

\def\lcm{\mathop{\rm lcm}\nolimits}

\def\Cyc{\mathop{\rm Cyc}}
\def\Bal{\mathop{\rm Bal}}

\def\lcc{\mathrel{\prec_{cc}}}

\def\eqcc{\mathrel{\approx_{cc}}}
\def\lecc{\mathrel{\preccurlyeq_{cc}}}

\def\leh{\mathrel{\preccurlyeq_{h}}}
\def\lh{\mathrel{\prec_{h}}}

\def\leTT{\mathrel{\preccurlyeq}}
\def\geTT{\mathrel{\succcurlyeq}}
\def\lTT{\mathrel{\prec}}
\def\gTT{\mathrel{\succ}}
\def\eqTT{\mathrel{\approx}}
\def\eqh{\mathrel{\approx_h}}

\DeclareMathAccent{\myarrow}{\mathord}{letters}{"7E}
\def\orC{\myarrow C}
\def\orG{\myarrow G}
\def\orH{\myarrow H}
\def\Ge{{\mathbb G}}
\def\calG{{\cal G}}
\newcommand{\TT}[1][]{\mathrel{\xrightarrow{\ifx @#1@ TT\else TT_{#1} \fi}}}
\newcommand{\nTT}[1][]{\thickspace\not\negthickspace\xrightarrow{\ifx @#1@ TT\else TT_{#1} \fi}}

\def\hom{\mathrel{\xrightarrow{hom}}}
\def\nhom{\mathrel{\thickspace\thickspace\not\negthickspace\negthickspace\xrightarrow{hom}}}
\def\floor #1{\left\lfloor #1 \right\rfloor}

\def\edge{\overrightarrow{K}_2}

\newenvironment{proofof}{\par\medskip\proofof}
  {\unskip\hfill$\Box$\par\medskip}
\def\proofof #1{\noindent{\bf Proof of #1:}}
\begin{document}

\def\makeheadbox{}
\title{Tension continuous maps---their structure and applications}
\author{
Jaroslav Ne\v set\v ril 
\thanks{Partially supported by COMBSTRU, by UPC, and by ICREA, Barcelona.}
 \and
Robert \v{S}\'amal
\thanks{Partially supported by UPC, Barcelona.}
}

\institute{Institute for Theoretical Computer Science (ITI)
  \thanks{Institute for Theoretical Computer Science is supported as
    projects LN00A056 and 1M0021620808 by Ministry of Education of the Czech Republic.}
 \\
 Charles University, Malostransk\'e n\'a\-m\v{e}s\-t\'{\i}~25, 
 118$\,$00 Prague, Czech Republic.
 \\
 \email {\{nesetril,samal\}@kam.mff.cuni.cz} }
\date{}
\maketitle
\begin{abstract}
We consider mappings between edge sets of graphs that lift
tensions to tensions. Such mappings are called tension-continuous
mappings (shortly $TT$~mappings). Existence of a $TT$~mapping
induces a (quasi)order on the class of graphs, which seems to 
be an essential extension of the homomorphism order
(studied extensively, see~\cite{HN}). In this paper we study the
relationship of the homomorphism and $TT$ orders.
We stress the similarities and the differences in both
deterministic and random setting.
Particularly, we prove that $TT$ order is dense and universal
and we solve a problem of M. DeVos et al. (\cite{DNR}).
\end{abstract}


\leftline{{\bfseries MSC\enspace}\enspace 05C15, 05C25, 05C38}

\section{Introduction} \label{sec:intro}



In this paper we study mappings between edge sets of graphs
that lift a tension to a tension. To motivate this we consider 
an important special case first. Let $G = (V,E)$ and $G'=(V',E')$ be
undirected graphs. A mapping $f : E \to E'$ is said to be 
\emph{cut-continuous} if for every cut $C$ in $G'$ the set $f^{-1}(C)$
forms a cut in~$G$.
(Here \emph{cut} means an edge cut, that is a set of all edges
between $X$ and $V\setminus X$ for some set $X \subseteq V(G)$.)
This condition is in particular satisfied when $f$ is induced by 
a homomorphism (see Lemma~\ref{homo}). 
However, there are other examples; in fact the main topic of this
paper is to understand to what extent these two notions coincide.

\begin{figure}[h]
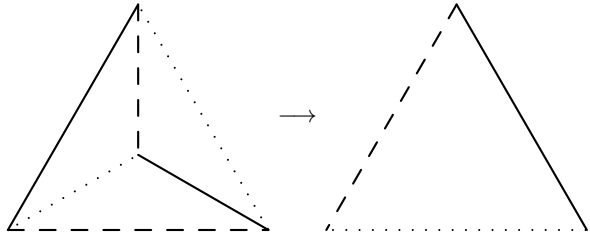

$$\vcenter{\hbox{\epsfbox{TT.4}}} 
\longrightarrow 
\vcenter{\hbox{\epsfbox{TT.3}}}$$
\caption{A cut-continuous mapping $K_4 \to K_3$ that is not induced by 
a mapping of vertices.}
\label{fig:K4}
\end{figure}

As a small example consider the 1-factorization of~$K_4$ 
(see Figure~\ref{fig:K4}). This constitutes a cut-continuous mapping 
$K_4 \to K_3$. On the other hand, an inclusion $K_3 \subseteq K_4$
is a homomorphism, hence induces a cut-continuous mapping. 
However, this example is an isolated one, as in Corollary~\ref{homotens}
we show that there are no other cut-continuous equivalent complete graphs.

Several other examples demonstrate that the existence of
cut-continuous mapping is not a very restrictive relation.
For example the well-known graphs depicted in Figure~\ref{fig:examples}
(on page~\pageref{fig:examples}) are all equivalent with respect to
cut-continuous mapping (see Theorem~\ref{examples}).
Another set of examples of cut-continuous mappings
we get by considering a pair of trees and any mapping
between their edge sets. 

This may indicate that the cut-continuous
mappings are abundant and induce a very weak graph comparison. 
Indeed, we give some further examples of pairs of graphs that
cannot be distinguished by cut-continuous mappings (although
they are not equivalent with respect to homomorphisms).
However, despite of all this evidence we prove that such examples
are rare, in the sense of random graphs.

Let a graph~$G$ be called \emph{homotens} if for every graph~$H$ every
cut-continuous mapping from~$G$ to~$H$ is induced by a homomorphism. In
Section~\ref{sec:random} we prove Corollary~\ref{homotens}, which in
particular implies the following theorem.

\begin {theorem}   \label{wbe}
A random graph is homotens with probability $1-o(1)$
(as size of the graph grows to infinity).
\end {theorem}

This result suggests to follow the now-standard approach to
homomorphisms (see e.g.~\cite{HN}) to investigate cut-continuous
mappings in the context of corresponding quasiorder~$\lecc$ and
strict partial order~$\lcc$. These are defined by
$$
   G \lecc G'  \hbox{ iff there is a cut-continuous mapping from~$G$ to~$G'$.}
$$
Also we let $G \eqcc G'$ denote $G \lecc G'$ and $G' \lecc G$.
The next theorem is an (important) special case
of Corollary~\ref{density} from Section~\ref{sec:density}.
(The density for homomorphic case was proven in~\cite{Wel}.)

\begin {theorem}[Density]   \label{jmeno}
For every pair of graphs $G_1 \lcc G_2$ (with the unique
exception $G_1 \eqcc K_1$, $G_2 \eqcc K_2$) there is a graph~$G$ such that
$$
    G_1 \lcc G \lcc G_2 \,.
$$
In other words, the order $\lcc$ is dense (if we do not consider
edgeless graphs).
\end {theorem}

Denote by $(Graphs,hom)$ ($(Graphs,cc)$, respectively)
the category of all finite graphs and all their homomorphisms
(all their cut-continuous mappings, respectively).
In Section~\ref{sec:universality} we prove Theorem~\ref{faithful}
that may be shortly expressed as follows.

\begin {theorem}   \label{embedding}
There is an embedding of $(Graphs,hom)$ into $(Graphs,cc)$.
\end {theorem}

\begin {corollary}[universality]   \label{universality}
Every countable partial order may be represented by (finite) graphs
with relation~$\lcc$.
\end {corollary}

The cut-continuous mappings were in the present context 
(that is as an important special case of tension-continuous mappings)
defined in~\cite{DNR}.
The motivation comes from Jaeger approach (\cite{Jaeger})
to classical conjectures (such as Berge-Fulkerson conjecture, 
Cycle Double Cover conjecture, Tutte's 5-flow conjecture).  
Let us remark that special cases of cut-continuous mappings were
(implicitly) studied earlier: 
\begin {itemize}
  \item Whitney classical theorem (\cite{Whit1}, \cite{Whit2}) 
    can be restated in our language: For 3-connected graphs $G$~and $G'$, 
    any bijection $f: E(G) \to E(G')$ such that both $f$~and $f^{-1}$ are
    cut-continuous 
  \footnote{or, equivalently, cycle-continuous---that is a preimage of a cycle is a cycle}
    is induced by
    an isomorphism. (A characterization for non-3-connected graphs is
    given as well.)
  \item  Kelmans (\cite{Kelmans}) generalized Whitney's theorem by
    introduction of cocircuit semi-isomorphisms
    of graphs. This is equivalent to our definition, although the
    notion is only defined when the mapping is a bijection.
  \item Linial, Meshulam, and Tarsi (\cite{LMT}) define cyclic (and 
    orientable cyclic) mappings. These are closely related to
    our definition of cut-continuous (and $\zet$-tension-continuous)
    mappings.
\end {itemize}
Our context is closest to that of~\cite{DNR}.
In~\cite{LMT}, only bijective mappings are considered, they serve as
a mean to define a variant of chromatic number ($\chi_{TT}$ of
Section~\ref{sec:codes}). We study non-bijective mappings too. This
(perhaps more natural) approach enables us to pursue the connections
between cut-continuous mappings and homomorphisms and
to study the properties of cut-continuous mappings in a broader view.
Other papers (\cite{Dvorak}, \cite{Reza}) will be mentioned at the relevant
place. For further research on this topic see~\cite{NS} and~\cite{RSthesis}.

The paper is organized as follows: In Section~\ref{sec:defce} we define 
group-valued tension-continuous mappings and prove their basic properties and
relevance to graph homomorphisms.  
We also briefly mention other types of XY-continuous mappings.
In Section~\ref{sec:density} we prove Density Theorem. 
This (perhaps surprisingly) relies on a new structural Ramsey-type theorem
(Theorem~\ref{Ramseybal}), which in turn leads to a solution of a problem
of~\cite{DNR}.  In Section~\ref{sec:random} we deal with random graphs and
as a consequence we are able to prove results analogous to the
homomorphism case (compare Theorem~\ref{homorigid} and
Corollary~\ref{TTrigid}). The properties of random graphs motivate
Section~\ref{sec:universality} where we prove the existence (by an
explicit construction) of rigid graphs (with respect to cut-continuous
mappings) and prove Theorem~\ref{embedding}.

Section~\ref{sec:group} is algebraic, we study the influence of a group~$M$
on the existence of $M$-tension-continuous mapping. The 
direct analogy of the Tutte result (dependence of the $M$-nowhere zero flow
only on~$|M|$) does not hold for $M$-tension-continuous mappings. Yet
we completely characterize the influence of the group in terms
of its algebraical structure (Theorem~\ref{numberofTT}).

In Section~\ref{sec:misc} we add several remarks and open 
problems. Particulary, we characterize 
(as a consequence of our approach) the complexity (and its 
dichotomy) of decision on the existence of a cut-continuous mapping.
Also, one has perhaps a surprising result that cut-continuous
mappings have no finite dualities (in the sense of~\cite{NT}).

\section{Definition \& Basic Properties}
\label{sec:defce}

\subsection{Basic notions}

We refer to \cite{Diestel}, \cite{HN} for basic notions on graphs and
their homomorphisms. 

By a graph we mean a finite\footnote{Although most of the results
  apply to infinite graphs too. The only place where infinite
  graphs appear is Corollary~\ref{homotens}, where we prove that
  tension-continuous mappings agree with homomorphisms on the
  Rado graph.}
directed graph, we write $uv$ (or sometimes $(u,v)$ for an edge 
from~$u$ to~$v$. (Occasionally we will speak of undirected graphs
too.) A \emph{circuit} in a graph is a connected subgraph in which
each vertex is adjacent to two edges. For a circuit~$C$, let
$C^+$ and~$C^-$ be the sets of edges oriented in either direction.
We will say that $(C^+, C^-)$ is a \emph{splitting} of edges of~$C$.

A \emph{cycle} is an edge-disjoint union of circuits. Given a graph~$G$
and a set~$X$ of its vertices, we let~$[X, \bar X]$ denote the set of all
edges with one end in~$X$ and the other in~$V(G) \setminus X$; we call
each such edge set a \emph{cut} in~$G$. Let $M$ be an abelian group.
We say that a function $\phi: E(G) \to M$ is an \emph{$M$-flow on~$G$}
if for every vertex
$v \in V(G)$  
$$
    \sum_{\hbox{$e$ enters $v$}} \phi(e) = 
    \sum_{\hbox{$e$ leaves $v$}} \phi(e) \,.
$$
Similarly, a function $\tau : E(G) \to M$ is an \emph{$M$-tension
on~$G$} if for every circuit~$C$ in~$G$ (with $(C^+,C^-)$ being
the splitting of its edges) we have
$$
  \sum_{e \in C^+} \tau(e) = \sum_{e \in C^-} \tau(e) \,.
$$

Note that $M$-tensions on a graph~$G$ form a vector space over~$M$
(if $M$ is a field, otherwise they only form a vector space
over~$\zet_2$), of dimension $|V(G)|-k(G)$, where $k(G)$ denotes
the number of components of~$G$. This vector space will be called
the $M$-tension space of~$G$; it is generated by \emph{elementary
$M$-tensions}, that is tensions that have some cut as their support.
(Elementary tensions are also called cut-tensions.)
Formally, for a
cut~$[X,\bar X]$ and $a \in M$ define
$$
   \phi^a_{[X,\bar X]} = \begin{cases}
     a  \quad \hbox{if $u \in X$ and $v \notin X$} \\
     -a  \quad \hbox{if $u \notin X$ and $v \in X$} \\
     0 \quad \hbox{otherwise}\,.\\
   \end{cases}
$$
Remark, that every $M$-tension is of form $\delta p$, 
where $p:V(G) \to M$ is any mapping and 
$(\delta p) (uv) = p(v) - p(u)$ (tension is a difference
of a potential).

Similarly, $M$-flows on~$G$ form a vector space of dimension
$|E(G)| - |V(G)| + k(G)$; it is generated by \emph{elementary
flows} (those with a circuit as a support) and it is orthogonal
to the $M$-tension space.

The above are the most basic notions in algebraic graph theory. For a
more thorough introduction to the subject see~\cite{Diestel}; we only mention
two most basic observations. 

A cycle can be characterized as a support of a $\zet_2$-flow and 
a cut as a support of a $\zet_2$-tension.
If $G$~is a plane graph then each cycle in~$G$ corresponds to a cut
in its dual~$G^*$; each flow on~$G$ corresponds to a tension on~$G^*$.

\subsection{Definitions}

The following is the principal notion of the paper. 
Let $G$, $G'$ be graphs and let $f: E(G) \to E(G')$ be a mapping
between their edge sets. 
We say $f$~is an \emph{$M$-tension-continuous mapping}
(shortly $TT_M$~mapping) if for every $M$-tension~$\tau$ on~$G'$, 
$\tau f$~is an $M$-tension on~$G$. The scheme below
illustrates this definition. It also shows that $f$ ``lifts
tensions to tensions'', explaining the term $TT$~mapping. 
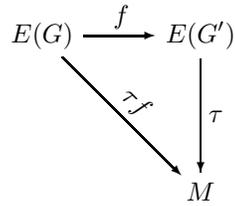
\begin{diagram}
      E(G) & \rTo^{f}       & E(G')     \\
           & \rdTo^{\tau f} & \dTo>\tau \\
           &                & M         \\
\end{diagram}
We write $f: G \TT[M] H$ if $f$ is a $TT_M$~mapping from~$G$
to~$H$ (or, more precisely, from~$E(G)$ to~$E(H)$).
In the important case $M = \zet_n$ we write $TT_n$ instead
of~$TT_{\zet_n}$. When $M$ is clear from the context, or when
we do not want to specify~$M$ we speak just of $TT$~mapping.

Let us only mention three related types of mappings:
$FF$ (lifts flows to flows), $FT$ (lifts tensions to flows), 
and $TF$ (lifts flows to tensions). 
In \cite{DNR} and \cite{RSthesis} these mappings are studied in more detail, 
in particular their connections to several classical 
conjectures (Cycle Double Cover conjecture, Tutte's 5-flow conjecture, and
Berge-Fulkerson matching conjecture) are explained. 

Of course if $M = \zet_2$ then the orientation of edges does not
matter. Hence, if $G$, $H$ are undirected graphs and 
$f: E(G) \to E(H)$ any mapping, we say that $f$ is 
\emph{$\zet_2$-tension-continuous} ($TT_2$) if for some
(equivalently, for every) orientation $\orG$ of~$G$ and $\orH$ of~$H$, 
$f$ is $TT_2$ mapping from~$\orG$ to~$\orH$. As cuts correspond
to $\zet_2$-tensions, with this provision $TT_2$ mappings of
undirected graphs are exactly cut-continuous mappings of
Section~\ref{sec:intro}.

Recall that $h: V(G) \to V(G')$ is called a \emph{homomorphism} if
for any $uv \in E(G)$ we have $f(u) f(v) \in E(G')$. We shortly
write $h: G \hom G'$ and define a quasiorder $\leh$ on the
class of all graphs by
$$
   G \leh G' \iff \hbox {there is a homomorphism $h: G \hom G'$.} 
$$
Homomorphisms generalize colorings: a $k$-coloring is exactly a
homomorphism $G \hom K_k$, hence $\chi(G) \le k$ iff $G \leh K_k$.
For an introduction to the theory of homomorphisms see \cite{HN}.

Motivated by the homomorphism order, we define for an abelian
group~$M$ an order $\leTT_M$ by
$$
   G \leTT_M G' \iff \hbox {there is a mapping $f: G \TT[M] G'$.} 
$$
This is indeed a quasiorder, see Lemma~\ref{compose}.
In Section~\ref{sec:group} we give a complete description of the
influence of the group $M$ on the notion of $M$-tension-continuous
mapping and on the relation $\leTT_M$.
We write $G \eqTT_M H$ iff $G \leTT_M H$ and $G \geTT_M H$, 
and similarly for~$G \eqh H$.
Ocassionally, we also write $G \TT[M] H$ instead of $G \leTT_M H$
and $G \hom H$ instead of $G \leh H$.

We define analogies of other notions used for study of homomorphisms:
a graph~$G$ is called \emph{$TT_M$-rigid} if there is no non-identical
mapping $G \TT[M] G$, graphs $G$, $H$ are called \emph{$TT_M$-incomparable}
if there is no mapping $G \TT[M] H$.

\subsection{Basic properties}

We start with an obvious yet key property of $TT_M$ mappings.

\begin {lemma}   \label{compose}
Let $f: G \TT[M]H$ and $g: H \TT[M] K$ be $TT_M$~mappings. Then the
composition $g \circ f$ is a $TT_M$~mapping.
\end {lemma}

\begin {lemma}   \label{TTsubgraph}
Let $f: G \TT[M] H$, let $H'$ be a subgraph of $H$, which contains
$f(e)$ for every $e \in G$. Then $f: G \to H'$ is $TT_M$ as well.
\end {lemma}

\begin {proof}
Take any $M$-tension $\tau'$ on~$H'$. Let $\tau' = \delta p'$
for $p': V(H') \to M$. If $V(H) = V(H')$ let $p = p'$, otherwise
extend~$p'$ arbitrarily to get~$p$. Now $\tau=\delta p$ is
an $M$-tension on~$H$ that agrees with~$\tau'$ on~$V(H')$. 
Hence $f \tau' = f \tau$, consequently $f \tau'$ is an~$M$-tension.
\qed
\end {proof}

\begin {corollary}   \label{factorization}
Let $f : G \TT[M] H$. Then there is a graph~$H'$ and 
$TT_M$ mappings $f_1: G \TT[M] H'$, $f_2: H' \TT[M] G$
such that $f_1$~is surjective and $f_2$~injective.
\end {corollary}

If $C$ is a circuit with a splitting~$(C^+, C^-)$, we say that 
$C$~is \emph{$M$-balanced} if for each $m\in M$ we have 
$(|C^+|-|C^-|)\cdot m = 0$. Otherwise, we say $C$~is \emph{$M$-unbalanced}.
We let $g_{\sss M}(G)$ denote the length of the shortest $M$-unbalanced
circuit in~$G$, or $\infty$ if there is none. 
For the particular case $M = \zet_2$ we can see that a circuit is $M$-balanced
if it is even, hence $g_{\sss \zet_2}(G)$ is the odd-girth of~$G$.
Easily, $G \TT[M] \edge$ iff any constant mapping $E(G) \to M$ is an
$M$-tension. This clearly happens precisely when all circuits in~$G$
are $M$-balanced, equivalently, if $g_{\sss M}(G) = \infty$.
As a consequence of this, the function $g_{\sss M}$ provides us with an
invariant for $TT_M$~mappings, as shown in the next two lemmas.

\begin {lemma}   \label{BalancedCircuit}
Let $M$ be an abelian group, let $G$, $H$ be graphs, let
$f : G \TT[M] H$. 
If $C$ is an $M$-unbalanced circuit in~$G$
then $f(C)$ contains an $M$-unbalanced circuit.
\end {lemma}

\begin {proof}
The inclusion homomorphism $C \to G$ induces a $TT_M$ mapping
$C \TT[M] H$. By Lemma~\ref{TTsubgraph} we get a mapping
$C \TT[M] f(C)$. If all circuits in~$f(C)$ are $M$-balanced, then
$f(C) \TT[M] \edge$ and, by composition we have $C \TT[M] \edge$.
This contradicts the fact that $C$ is $M$-unbalanced.
\qed
\end {proof}

\begin {lemma}   \label{invariant}
Let $G \leTT_M H$. Then $g_{\sss M}(G) \ge g_{\sss M}(H)$. 
\end {lemma}

\begin {proof}
If $g_{\sss M}(G)=\infty$, the conclusion holds. Otherwise, let $C$
be an $M$-unbalanced circuit of length~$g_{\sss M}(G)$ in~$G$. 
By Lemma~\ref{BalancedCircuit}, $f(C)$ contains an $M$-unbalanced
circuit. It is of size at least~$g_{\sss M}(H)$ and at most~$g_{\sss M}(G)$.
\qed
\end {proof}

For a homomorphism $h: G \to G'$ we write $h^\sharp$ for the
\emph{induced mapping on edges}, that is $h^\sharp (uv) = h(u)h(v)$. 
The following easy lemma is the starting point of our investigation.

\begin {lemma}   \label{homo}
Let $G$, $H$ be graphs, $M$ abelian group.
For every homomorphism $f : G \hom H$ the induced mapping
$f^\sharp : E(G) \to E(H)$ is $M$-tension-continuous
(in particular cut-continuous).
Hence, from $G \leh H$ follows $G \leTT_M H$.

If $f: V(G) \to V(H)$ is an antihomomorphism (that is, 
it reverses every edge), $f^\sharp$ is $M$-tension-continuous, too.
\end {lemma}

\begin{proof}
Let $f: G \to H$ be a homomorphism, $\phi: V(H) \to M$ a tension. We may
assume that $\phi$ is a cut-tension corresponding to the cut 
$[X, V\setminus X]$. Then the cut $[f^{-1}(X),f^{-1}(V\setminus X)]$
determines precisely the tension $\phi\circ f$.
\qed
\end{proof}

The main theme of this paper is to find similarities and differences
between orders $\leh$ and $\leTT_M$. 
In particular we are interested in when the converse to Lemma~\ref{homo}
holds: 

\begin{problem} \label{homoTT}
Let $f : E(G) \to E(H)$. Find suitable conditions for $f$, $G$, $H$ 
that will guarantee that whenever $f$ is $TT_M$, then it
is induced by a homomorphism (or an antihomomorphism); 
i.e.\ that there is a homomorphism (or an antihomomorphism) 
$g : V(G) \to V(H)$ such that $f = g^\sharp$.
\end{problem}

Applying no further conditions this does not hold, see examples in
the first section, Theorem~\ref{examples}, and
Theorem~\ref{largeexamples}. We start with a result
that provides a condition on~$H$ (cf.~\cite{DNR}).
If $M$ is any group and $B \subseteq M \setminus \{0\}$
any set, then we define
$$
   \Cay(M,B) = (M, \{uv, v-u\in B\}).
$$

\begin {lemma}   \label{Hcayley}
Let $M$ be a group, $B \subseteq M^n \setminus \{(0,\dots, 0)\}$.
Denote $H = \Cay(M^n, B)$. Then for every graph~$G$
$$
    G \TT[M] H    \iff    G \hom H \,,
$$
in fact every $TT_M$ mapping is induced by a homomorphism
or by an antihomomorphism.
\end {lemma}


Another partial answer to Problem~\ref{homoTT} is to put 
some restriction on~$G$. This seems more fruitful as
the necessary restriction is rather weak. 
We will say that~$G$ is \emph{$M$-homotens} if 
for any graph~$H$ any $TT_M$ (that is $M$-tension continuous) mapping
$G \TT[M] H$ is induced by a homomorphism (or an antihomomorphism). 
Note that if $M = \zet_2^n$, all $2^{|E(G)|}$ orientations of a
graph~$G$ are $TT_M$-equivalent. Thus, for such $M$ it makes sense
to investigate $M$-homotens undirected (instead of directed)
graphs: We say an undirected graph~$G$ is $\zet_2$-homotens
if for any undirected graph~$H$ any $TT_2$ mapping $G \TT[2] H$
is induced by a homomorphism of the undirected graphs.

As we deal mostly with the case $M = \zet_2$, we call $\zet_2$-homotens
graphs shortly homotens. In Section~\ref{sec:random} we prove a perhaps
surprising fact that most of the graphs are homotens.

Yet another partial answer to Problem~\ref{homoTT} 
is provided by non-trivial theorem (proved in~\cite{DNR})
that studies mappings, which are defined more
restrictively than $TT$. In other words, we put restrictions
on~$f$ this time.

A mapping $f: E(G)\to E(G')$ is \emph{$\zet$-cut-tension-continuous}
iff for every cut-tension $\phi: E(G') \to \zet$ the mapping
$\phi \circ f$ is a cut-tension $E(G) \to \zet$.

\begin{theorem}[\cite{DNR}] \label{cut} 
Any $\zet$-cut-tension-continuous mapping
is induced by a homomorphism or by an antihomomorphism.
\end{theorem}

Before proceeding any further we present an alternative definition 
of tension-continuous mappings (which is proved in~\cite{DNR}).
For mappings $f: E(G) \to E(H)$ and $\phi: E(G) \to M$ we let $\phi_f$
denote the \emph{algebraical image of $\phi$}: that is we define a mapping
$\phi_f:E(H) \to M$ by 
$$
  \phi_f(e') = \sum_{e\in E(G); f(e) = e'} \phi(e) \,.
$$

\begin {lemma}   \label{altdef}
Let $f: E(G) \to E(H)$ be a mapping. Then $f$ is
$M$-tension-continuous if and only if for every
$M$-flow $\phi$ on~$G$, its algebraical image $\phi_f$
is an~$M$-flow. 

We formulate this explicitly for $M=\zet_2$: Mapping $f$ is cut-continuous 
if and only if for every cycle~$C$ in~$G$, the set of edges
of~$H$, to which an odd number of edges of~$C$ maps, is a cycle.
\end {lemma}

The following interesting construction provides a completely 
different connection between homomorphisms and tension-continuous
mappings. Given an (undirected) graph $G = (V,E)$ write $\Delta(G)$
for the graph $({\cal P}(V),\Delta(E))$, where $AB \in \Delta(E)$ iff
$A \Delta B \in E$ (here ${\cal P}(V)$ denotes the set of all
subsets of $V$ and $A \Delta B$ the symmetric difference of sets
$A$ and $B$).

\begin{theorem} \label{cube} 
Let $G$, $H$ be undirected graphs. Then $G \leTT_{\zet_2} H$ 
iff $G \leh \Delta(H)$.
\end{theorem}

We could formulate an analogous construction and result for groups
$M \ne \zet_2$; the role of~$\Delta(H)$ would be played by some
Cayley graph on the group $M^n$ for appropriate~$n$; 
for finite~$M$, this Cayley graph is finite.
Theorem~\ref{cube} is proved in~\cite{DNR} and (for $H = K_n$)
in~\cite{LMT}. As we shall make use of it we prove it here for the sake of
completeness.

\begin{proof}
Given $g: G \TT[2] H$ we construct $f: G \hom \Delta(H)$ as follows.
First choose $v_0 \in V(G)$ and let $f(v_0) = \emptyset$. Then
whenever $uv$ is an edge with $f(v)$ defined, 
we let $f(u) = f(v) \Delta g(uv)$. Using Lemma~\ref{altdef} it can
be easily verified that the construction is consistent, clearly
it defines a homomorphism.
For the backward implication we just define $g(uv) = f(u)\Delta f(v)$
and apply Lemma~\ref{altdef}.
\qed
\end{proof}

These two results may indicate that quasiorders $\leTT_M$ and
$\leh$ are closely related. Before pursuing the similarities, 
we stress some of the differences.

\begin{figure}
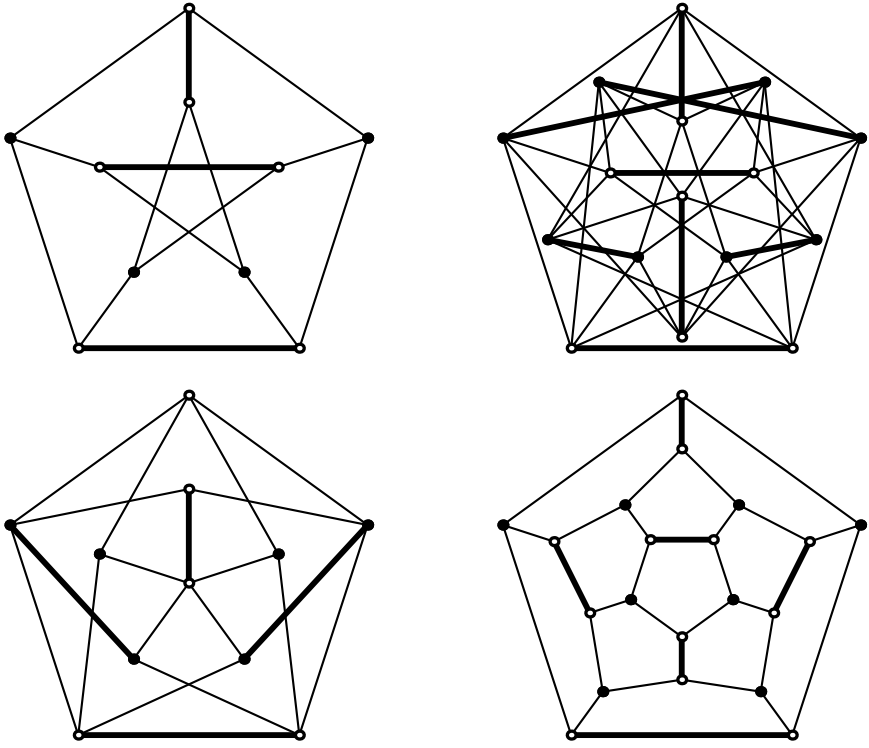

\hbox to \hsize{\epsfbox{Graphs.0}  \hfill  \epsfbox{Graphs.1}}

\bigskip
\hbox to \hsize{\epsfbox{Graphs.2}  \hfill  \epsfbox{Graphs.3}}

\caption{Examples of graphs that are $TT$-equivalent to $C_5$.
One color class is drawn in bold, the other four are obtained
by rotation.}
\label{fig:examples}
\end{figure}

\begin {theorem}   \label{examples}
Let $P$ be the Petersen graph, $Cl$ the Clebsch graph, $Gr$ the
Gr\" otsch graph, $D$ the dodecahedron (see Figure~\ref{fig:examples}).
Then $P \eqTT_{\zet_2} Cl \eqTT_{\zet_2} Gr \eqTT_{\zet_2} 
D \eqTT_{\zet_2} C_5$. On the other hand, in the homomorphism
order no two of these graphs are equivalent.
\end {theorem}

\begin{proof}
We have $C_5 \subset P \subset Cl$, $C_5 \subset D$, and $C_5 \subset Gr$.
As inclusion is a homomorphism and hence it induces a $TT$~mapping, we
only need to provide mappings $Cl \TT C_5$, $D \TT C_5$, and $Gr \TT C_5$. 
In Figure~\ref{fig:examples}, we emphasize some edges in each graph.
Let $G$ be the considered graph and $A \subseteq E(G)$ the set
of bold edges. Put $A_1 = A$ and let $A_2$, $A_3$, $A_4$,
$A_5$ denote the sets obtained from~$A$ by rotation, so that the
sets $A_i$ partition $E(G)$. Define a mapping 
$E(G) \to E(C_5) = \{e_1, \dots, e_5\}$ 
by sending all edges in~$A_i$ to $e_i$.

Note that 4-edge subgraphs of $C_5$ generate its
$\zet_2$-tension space. Hence it is enough to verify that after
deleting any color class we are left with a cut. Due to symmetry 
we only need to check that $E(G)\setminus A$ is a cut in~$G$. 
This is straightforward to verify, the corresponding bipartition
of vertices is depicted in Figure~\ref{fig:examples}.
\qed
\end{proof}

Graphs $\Delta(K_n)$ will be further studied in Section~\ref{sec:misc}.
Here we only illustrate Theorem~\ref{cube} by a particular choice $H=C_5$. 
Graph $\Delta(H)$ consists of two components, each of
which is isomorphic to the Clebsch graph~$Cl$. 
Hence, $G \TT C_5$ is equivalent to $G \hom Cl$. 
This reproves part of Theorem~\ref{examples} but, more importantly,
this observation implicitly appeared in~\cite{Reza},
where a theorem of~\cite{Gue} was used to prove the following result.

\begin {theorem}   \label{reza-thm}
Any planar triangle-free graph admits a homomorphism
to the Clebsch graph.
\end {theorem}

The next theorem gives an infinite class of graphs where homomorphisms
and $TT$~mappings differ. In particular it implies that for every $n$
there are $n$-connected graphs that are not homotens.

\begin {theorem}   \label{largeexamples}
Let $n$ be odd. Denote $G_n$ one of the (two isomorphic)
components of~$\Delta(K_n)$.
Graphs $K_n$ and $G_n$ are $TT_2$-equivalent and both 
are $(n-1)$-connected.
Finally, $G_n \nhom K_n$ for $n=2^k-1$.
\end {theorem}

\begin {proof}
Using Theorem~\ref{cube} for $G = H = K_n$ we get 
$K_n \hom \Delta(K_n)$. From connectivity of~$K_n$ and
from Lemma~\ref{homo} it follows $K_n \TT[2] G_n$.
Using Theorem~\ref{cube} for $G = H = \Delta(K_n)$
we get $\Delta(K_n) \TT[2] K_n$, hence also $G_n \TT[2] K_n$.

Graph $K_n$ is $(n-1)$-connected. 
Easily $\Delta(K_n) = Q_n^{(2)}$, where $Q_n$~is the $n$-dimensional
hypercube and $Q_n^{(2)}$ means that we are connecting by an edge the
vertices at distance two in the hypercube. It is well-known
and straightforward to verify that $Q_n$ is $(n-1)$-connected.
The vertices with odd (even) number of 1's among their coordinates
form the two components of~$Q_n^{(2)}$; for an odd~$n$ these two components 
are isomorphic by a mapping $\vec x \mapsto (1,1,\dots, 1) - \vec x$.
Observe that if we take a path in~$Q_n$ and leave every second
vertex out, we obtain a path in~$Q_n^{(2)}$. So $Q_n^{(2)}$ is
$(n-1)$-connected since $Q_n$ is.

For the last part of the theorem, it follows 
from the remarks in the Section~\ref{sec:codes} 
that $\chi(G_n) = n+1$ for $n = 2^k - 1$.
\qed
\end {proof}

\section {Density}
\label{sec:density}
\subsection {A Ramsey-type theorem for locally balanced graphs}
\label{sec:Ramsey}

In this subsection we deal with undirected graphs only.
We prove a Ramsey-type theorem that will be used in
Section~\ref{sec:SIL} as a tool to study $\lTT_M$ (on directed
graphs).

An \emph{ordered graph} is an undirected graph with a fixed linear
ordering of its vertices. The ordering will be denoted by $<$, an
ordered graph by~$(G,<)$, or shortly by~$G$.
We say that two ordered graphs are \emph{isomorphic},
if the (unique) order-preserving bijection is a graph isomorphism.
An ordered graph $(G,<)$ is said to be a \emph{subgraph} of~$(H,<')$, if
$G$~is a subgraph of~$H$, and the two orderings coincide on $V(G)$.

A circuit $C = v_1, \dots, v_l$ in an ordered graph is \emph{balanced}
iff 
$$
  |\{ i ; v_i < v_{(i \bmod l) + 1}\}| = 
  |\{ i ; v_i > v_{(i \bmod l) + 1}\}| \,.
$$
This can be reformulated using the notion preceding 
Lemma~\ref{BalancedCircuit}. Let $\orG$ be a directed graph
with $V(\orG) = V(G)$ and 
$E(\orG) = \{ (u,v); \mbox{$uv \in E(G)$ and $u < v$}\}$.
(We can say that all edges are oriented ``up''.)
Then a circuit in~$G$ is balanced iff the corresponding circuit
in~$\orG$ is $\zet$-balanced. Note that a circuit in~$\orG$
is $\zet_2$-balanced iff its length is even.

Denote by $\Cyc_p$ the set of all ordered graphs that contain no 
odd circuit of length at most $p$.
Denote by $\Bal_p$ the set of all ordered graphs that contain no 
unbalanced circuit of length at most $p$. 

Ne\v set\v ril and~R\"odl (\cite{NR}) proved the following Ramsey-type
theorem.

\begin{theorem}   \label{Ramsey}
Let $k$,~$p$ be positive integers. For any ordered graph 
$(G,<) \in \Cyc_p$ there is an ordered graph $(H,<) \in \Cyc_p$ with the
``Ramsey property'': for every coloring of~$E(H)$ by $k$
colors there is a monochromatic subgraph ${(G',<)}$, isomorphic to~$(G,<)$.
\end {theorem}

We will need a version of this theorem for $\Bal_p$.
By the discussion above, this means that we consider $\zet$-balanced
(instead of $\zet_2$-balanced) circuits.

\begin{theorem}   \label{Ramseybal}
Let $r$,~$p$ be positive integers. For any ordered graph 
$(G,<) \in \Bal_p$ there is an ordered graph $(H,<) \in \Bal_p$ with
the ``Ramsey property'': for every edge coloring of $H$ by
$r$ colors there is a monochromatic subgraph ${(G',<)}$, isomorphic
to~$G$. This conclusion will be shortly written as $(H,<) \to (G,<)^2_r$.
\end {theorem}

\begin {proof}
The proof of Theorem~\ref{Ramseybal} uses a variant of the
amalgamation method (partite construction) due to the
first author and R\"odl (see e.g. \cite{NR-Ramsey}, \cite{Nes-Ramsey}), 
which has many applications in structural Ramsey theory. 

For the purpose of this proof we slightly generalize the notion 
of ordered graph. 
We work with graphs with a quasiordering $\le$ of its vertices; 
such graphs are called \emph{quasigraphs}, $\le$ is called the
standard ordering of~$G$. Alternatively, a quasigraph $(G, \le)$
is a graph $G = (V,E)$ with a partition 
$V_1 \cup V_2 \cup \cdots \cup V_a$ of $V$: each $V_i$ is a set
of mutually equivalent vertices of~$V$ and $V_1 < V_2 < \cdots < V_a$.
The number $a$ of equivalence classes of $\le$ will be fixed
throughout the whole proof. In this case we speak about
\emph{$a$-quasigraphs}. It will be always the case that every $V_i$
is an independent set of~$G$.

An embedding $f: (G, \le) \to (G', \le')$  is an embedding 
(i.e. an isomorphism onto an induced subgraph) $G \to G'$
which is moreover monotone with respect to the standard orderings
$\le$ and~$\le'$. Explicitly, such an embedding~$f$ is an embedding
of~$G$ to~$G'$ for which there exists an increasing mapping 
$\iota: \{ 1, 2, \dots, a\} \to \{ 1, 2, \dots, a'\}$ such that
$f(V_i) \subseteq V'_{\iota(i)}$ for $i = 1, \dots, a$. 
(Here $V'_1 < V'_2 < \cdots < V'_{a'}$ are equivalence classes
of the quasiorder $\le'$.) By identifying the equivalent vertices of a
quasigraph~$G$ we get a graph~$\bar G$ and a homomorphism $\pi:G \to \bar G$;
graph $\bar G$ is called the \emph{shadow} of~$G$, mapping $\pi$ is 
called \emph{shadow projection}.

We prove Theorem~\ref{Ramseybal} by induction on~$p$. 
The case $p=1$ is the Ramsey theorem for ordered graphs and so
we can use Theorem~\ref{Ramsey} for~$p=1$. In the induction step
($p \to p+1$) consider arbitrary ordered graph $(G, \le)$, 
let $G = (V,E)$, $|V| = n$, and $G \in \Bal_{p+1}$.
By the induction assumption there exists an ordered graph
$(K,\le) \in \Bal_p$ such that 
$$
  K \to (G)_r^2  \,.
$$
Let $V(K) = \{ x_1 < \cdots < x_a\}$
and $E(K) = \{ e_1, \dots, e_b\}$. In this situation 
we shall construct (by induction) $a$-quasigraphs 
$P^0$, $P^1$, \dots, $P^b$ (called usually ``pictures''). Then the
quasigraph~$P^b$ will be transformed to the desired ordered 
graph~$(H, \le) \in \Bal_p$ satisfying
$$
    (H,\le) \to (G, \le)_r^2 \,.
$$
We proceed as follows. Let $(P^0, \le^0) \in \Bal_{p+1}$ be
$a$-quasiordered graph for which for every induced 
subgraph~$G'$ of~$K$, such that $(G',\le)$ is isomorphic to $(G,\le)$,
there exists a subgraph~$G_0$ of~$P^0$ with the shadow~$G'$. 
Clearly $(P^0, \le^0)$ exists, as it
can be formed by a disjoint union of~$\binom an$ copies of~$G$ with an
appropriate quasiordering. 

In the induction step $k \to k+1$ ($k \ge 0$) let the picture
$(P^k, \le^k)$ be given. Write $P^k = (V^k, E^k)$
and let $V_1^k < V_2^k < \cdots < V_a^k$ be all equivalence classes
of~$\le^k$. Consider the edge $e_{k+1} = \{x_{i_{k+1}}, x_{j_{k+1}}\}$ of~$K$
($x_{i_{k+1}} < x_{j_{k+1}}$).
To simplify the notation, we will write $i=i_{k+1}$, $j=j_{k+1}$.
Let $B^k = (V^k_i \cup V^k_j, F^k)$ be the bipartite subgraph of~$P^k$
induced by the set $V^k_i \cup V^k_j$. We shall make use of the
following lemma.

\begin {lemma}   \label{bipartite:Ramsey}
For every bipartite graph~$B$ there exists a bipartite
graph~$B'$ such that
$$
        B' \to (B)_r^2 \,.
$$
(The embeddings of bipartite graphs map the upper part to the upper
part and the lower part to the lower part.)
\end {lemma}

This lemma is easy to prove and it is well-known, see
e.g.~\cite{Nes-Ramsey}.

Continuing our proof, let 
\begin{equation}\label{eq:bipartite}
  B'^k \to (B^k)_r^2
\end{equation}
be as in Lemma~\ref{bipartite:Ramsey} and put explicitly 
$B'^k = (V_i^{k+1} \cup V_j^{k+1}, F^{k+1})$. Let also 
${\cal B}_k$ be the set of all induced subgraphs of~$B'^k$, which are
isomorphic to~$B^k$. Now we are in the position to construct the
picture $(P^{k+1}, \le^{k+1})$.

We enlarge every copy of~$B^k$ to a copy of~$(P^{k}, \le^k)$ while
keeping the copies of~$P^k$ disjoint outside the set 
$V_i^{k+1} \cup V_j^{k+1}$. The quasiorder $\le^{k+1}$ is defined
from copies of quasiorder $\le^k$ by unifying the corresponding classes.
While this description perhaps suffices to many here is an explicit
definition of~$P^{k+1}$:

Put $P^{k+1} = (V^{k+1}, E^{k+1})$, where 
$V^{k+1} = V^k \times {\cal B}/{\sim}$. The equivalence $\sim$ 
is defined by 
$$
  (v, B) \sim (v', B') 
    \quad
    \Longleftrightarrow 
    \quad
    v = v' \in V_i^{k+1} \cup V_j^{k+1} 
    \quad \hbox{or $\quad v = v'$ and $B=B'$.}
$$
Denote by $[v,B]$ the equivalence class of $\sim$ containing
$(v,B)$. We define the edge set by 
putting $\{[v,B],[v',B']\} \in E^{k+1}$ if $\{v,v'\} \in E^k$
and $B = B'$. Define quasiorder $\le^{k+1}$ by putting
$$
   [v,B] \le^{k+1} [v',B'] \Longleftrightarrow v \le^k v' \,.
$$
It follows that $\le^{k+1}$ has $a$ equivalence classes
$V_1^{k+1} < \cdots < V_a^{k+1}$. (Note that this is
consistent with the notation of classes $V_i^{k+1}$, $V_j^{k+1}$
of~$B'^k$.) 

Continuing this way, we finally define the picture $(P^b, \le^b)$.
Put $H = P^b$ and let $\le$ be an arbitrary linear ordering that
extends the non-symmetric part of the quasiorder $\le^b$. We claim
that the graph $H$ has the desired properties. To verify this it
suffices to prove:
\begin {enumerate}[(i)]
  \item $(H,\le) \in \Bal_{p+1}$ and
  \item $(H,\le) \to (G,\le)_r^2$.
\end {enumerate}

The statement (i) will be implied by the following claim.

\begin{claim}
\begin {enumerate}
  \item $(P^0, \le^0) \in \Bal_{p+1}$. 
  \item If $(P^k, \le^k) \in \Bal_{p+1}$, then ${(P^{k+1},\le^{k+1})} \in \Bal_{p+1}$.
\end {enumerate}
\end{claim}


\begin{proof}[of Claim]
The first part follows from the construction.
In the second part, suppose that $P^{k+1}$~contains an unbalanced
circuit $C = u_1,u_2, \dots, u_l$ of length $l \le p+1$.
Let $\pi: V(P^{k+1}) \to V(K)$ be the projection, that is
for $u \in V_s^{k+1}$ we have $\pi(u) = x_s$. From the construction
it follows that $\pi$~is a homomorphism $P^{k+1} \hom K$, in other words
that $K$~is the shadow of~$P^{k+1}$.

Consider the closed walk $C_\pi = \pi(u_1), \pi(u_2), \dots, \pi(u_l)$
in~$K$.  As $C_\pi$ is unbalanced closed walk, it contains an
unbalanced circuit of length~$l' \le l$. Since $K \in \Bal_p$, we have 
$l'=l=p+1$, that is $\pi(u_1), \dots, \pi(u_{p+1})$ are all distinct.
Let $u_s = [v_s,B_s]$. If $B_1 = B_2 = \cdots = B_l$, that is the
whole~$C$ is contained in one copy of~$P^k$, we have a contradiction
as $P^k \in \Bal_{p+1}$. 

Now we use the construction of~$P^{k+1}$ as an amalgamation of
copies of~$P^k$: If $B_t \ne B_{t+1}$ (indices modulo~$l$),
then $\pi(u_t) \in \{x_{i_{k+1}}, x_{j_{k+1}}\}$. As the vertices 
$\pi(u_1), \dots, \pi(u_l)$ are pairwise distinct, this happens just
for two values of~$t$. Consequently, the whole~$C'$ is contained
in two copies of~$P^k$ and there are indices $\alpha$, $\beta$ such
that $\pi(u_\alpha) = x_{i_{k+1}}$ and $\pi(u_\beta) = x_{j_{k+1}}$.

The circuit~$C$ is a concatenation of $P'$ and~$P''$---two paths 
between $u_\alpha$ and $u_\beta$, each of them properly contained
in one copy of~$P^k$. No copy of~$P^k$ contains whole~$C$, therefore
both $P'$ and~$P''$ have at least two edges, hence 
at most $p-1$ edges. Let $\bar P'$, $\bar P''$ denote the shadows
of~$P'$ and~$P''$. Both 
$\bar P' \cup \{e_{k+1}\}$ and~$\bar P'' \cup \{e_{k+1}\}$
are closed walks in~$K$ containing at most~$p$ edges. As $K \in \Bal_p$, 
both of them are balanced, so $C$ is balanced as well, a
contradiction.
\qed
\end{proof}

\def\aa{{\cal A}}
We turn to the proof of statement (ii). We use a standard argument that
is the core of the amalgamation method. 
Let $E(H) = E(P^b) = \aa_1 \cup \dots \cup \aa_r$ be a fixed coloring.
We proceed by backwards induction $b \to b-1 \to \cdots$ and we prove
that there exists a quasisubgraph 
$P_0^k$ of~$P^b$ isomorphic to~$P^k$ such that for any $l>k$, any two edges 
of~$P_0^k$ with shadow $e_l$ get the same color.
This is easy to achieve using the Ramsey properties~(\ref{eq:bipartite}) of
graphs~$B'^k$. Finally, we obtain a copy $P_0^0$
of~$P^0$ in~$P^b$ such that the color of any of its edges 
depends only on its shadow (in~$K$). However $K \to (G)_r^2$
and as for any copy~$G'$ of~$G$ in~$K$ there exists a subgraph 
$G_0$ of~$\bar P^0$ such that its shadow is~$G'$ we get 
that there exists a monochromatic copy of~$G$ in~$P^b$. This
concludes the proof.
\qed
\end {proof}

\subsection{Density}
\label{sec:SIL}

In this section we prove the density of $TT_M$ order
(for every abelian group~$M$).  For this we first prove the ``Sparse
Incomparability Lemma'', Lemma~\ref{SIL} (analogous statement for
homomorphisms appears in~\cite{NR-rigid}). Although the proof follows
similar path as in the homomorphism case, some steps are considerably
harder; the main reason for this is the nonexistence of products
in the category of tension-continuous mappings. To overcome this
obstacle, we use the Ramsey-type theorem from the previous subsection.

\begin{figure}
\epsfbox{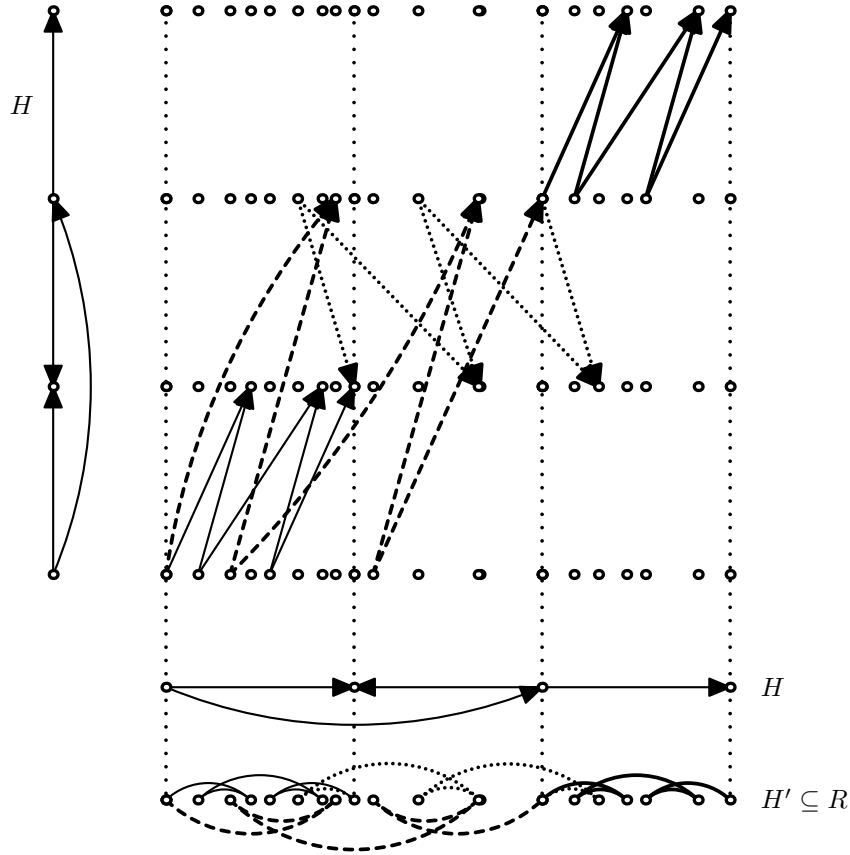}
\caption{An illustration of the proof of Lemma~\ref{SIL}
(here $s=5$). Only a part of the graph~$G' = H\times R$ is shown.}
\label{fig:SIL}
\end{figure}

\begin {lemma}   \label{SIL}
Let $M$ be an abelian group, let $l$, $t\ge 1$ be integers.
Let $G_1$, $G_2$,~\dots, $G_t$, $H$ be graphs such that
$H \nTT[M] G_i$ for every~$i$ and $H \nTT[M] \edge$.
Then there is a graph~$G'$ such that
\begin {enumerate}
\item $G' \lTT_M H$, moreover $G' \lh H$, 
\item all circuits in $G'$ shorter than~$l$ are $M$-balanced, and
\item $G' \nTT[M] G_i$ for every $i = 1, \dots, t$.
\end {enumerate}
\end {lemma}

\begin {proof}
Choose an odd integer $p$ larger than $\max \{ |E(H)|, l \}$. 
Pick any linear ordering of $V(H)$ to make~$H$ into an ordered
graph~$(H,<)$ and subdivide each edge to increase the girth.
More precisely, we replace every edge~$e$ of~$H$ by an oriented
path~$P(e)=e_1,e_2,\dots, e_p$; 
the ordering of $V(H)$ is extended to the new vertices so, that 
$e_j$ goes up iff $j$ is odd, see Figure~\ref{fig:SIL}.
When we do this for every edge of~$H$, we forget the orientation of
the edges and let~$(H',<)$ denote the resulting ordered graph. It
is $(H',<)\in \Bal_p$.

Put $r = \max_i |E(G_i)|^{|E(H)|}$. Using Theorem~\ref{Ramseybal} 
we find a graph~$(R,<)\Bal_p$ satisfying $(R,<) \to (H',<)^2_r$
As every circuit of~$(R,<)$ is balanced, it is also $M$-balanced.

We orient all edges of~$R$ up (that is towards the vertex larger
in~$<$), and set~$G'= H\times R$ 
(see Figure~\ref{fig:SIL}).
Formally, $V(G') = V(H) \times V(R)$, and for edges $e=uv$ of~$H$ and
$e'=u'v'$ of~$R$ we have an edge from~$(u,u')$ to~$(v,v')$
(this edge will be denoted by~$(e,e')$).

\def\Hbar{\bar H}
Now $G' \TT[M] R$ (as there is even a homomorphism---the projection),
so by Lemma~\ref{BalancedCircuit} there is no short $M$-unbalanced circuit
in~$G'$. This gives part~2 of the statement.
The other projection of~$G'$ gives $G' \TT[M] H$, and indeed
even $G' \hom H$.  To prove part~1, 
we need to exclude the case $H \TT[M] G'$. If such a mapping
exists, denote $\Hbar$ its image in~$G'$. It is easy to verify that
$H \TT[M] \Hbar$. As $s > |E(H)|$, there is no $M$-unbalanced
circuit in~$\Hbar$, hence $\Hbar \TT[M] \edge$. By composition
we get $H \TT[M] \edge$, a contradiction. Hence $G' \lTT_M H$, and
therefore $G' \lh H$ too. It remains to prove part~3.

For the contrary, suppose there is an index $i$ and a $TT_M$
mapping $f : G' \TT[M] G_i$. As $G' = H \times R$, this induces
a coloring~$c$ of edges of~$R$ by elements of~$E(G_i)^{E(H)}$
(where $c(e')$ sends $e$ to $f((e,e'))$). As we have chosen $R$ to be
a Ramsey graph for $H'$, there is a monochromatic copy of~$H'$ in~$R$.
To ease the notation we will suppose this copy is just~$H'$, let $g$
be the color of edges of~$H'$. We will show that $g$~is 
a $TT_M$~mapping $H \to G_i$, and this will be our desired contradiction.

\def\Cbar{\bar C}
\def\phibar{\bar \phi}
We will use Lemma~\ref{altdef}, hence for any flow $\phi:E(H)\to M$
we need to show that $\phi_g$ is a flow. 
Clearly it is enough to verify this for $\phi$ being an elementary
flow, as elementary flows generate the $M$-flow space on~$H$.
So let $C$~be a circuit in~$H$ that is the support of~$\phi$.
The corresponding circuit $\Cbar$ in $H \times R$ has edge set
$$
  E(\Cbar) = \bigcup \{ \{e\}\times P(e), e\in E(C) \} \,.
$$
Let $\phibar$ be the elementary flow on $H \times R$ corresponding
to~$\phi$. Explicitly, 
$$
   \phibar : (e,e_i) \mapsto 
    \begin{cases}
       \phi (e)   & \hbox{if $i$ is odd,} \\
      -\phi (e)   & \hbox{if $i$ is even.} \\
    \end{cases}
$$
As $H'$ is $g$-monochromatic, $f((e,e_i)) = g(e)$ for every $i$.
Consequently $\phi_g = \phibar_f$, so $\phi_g$ is a flow.  
\qed
\end {proof}

\begin {theorem}   \label{AntiExtInt}
Let $M$ be an abelian group, let $t \ge 0$ be an integer.
Let $G$, $H$ be graphs such that $G \lTT_M H$ and $H \nTT[M] \edge$. 
Let $G_1$, $G_2$, \dots, $G_t$ be pairwise
incomparable (in $\lTT_M$) graphs satisfying $G\lTT_M G_i \lTT_M H$ for
every $i$. Then there is a graph~$K$ such that
\begin {enumerate}
  \item $G \lTT_M K \lTT_M H$,
  \item $K \nTT G_i \nTT K$ for every $i = 1, \dots, t$.
\end {enumerate}
If in addition $G \hom H$ then we have even $G \lh K \lh H$.
\end {theorem}

\begin {proof}
Choose $l > \max \{ |E(H)|,  |E(G_i)|, i=1, \dots, t\}$.
We use Lemma~\ref{SIL} to get a graph~$G'$ such that 
$G' \nTT G_i$ and $G' \nTT G$; then put $K = G + G'$. 
Easily $G \leTT K \le H$ and $K \nTT G_i$, $K \nTT G$
(as $G'$ has this property). It remains to show $F \nTT K$
for $F \in \{ H, G_1, \dots, G_t\}$. Note that it is
not enough to show $F \nTT G$ and $F \nTT G_i$, we have
to proceed more carefully. 

So suppose we have an $TT_M$~mapping $f: F \TT G + G'$. Pick an edge 
$e_0 \in E(G)$, and define $g: E(F) \to E(G)$ as follows:
$$
  g(e) =
  \begin{cases}
    f(e) & \text{if $f(e) \in E(G)$} \\
    e_0  & \text{otherwise.} \\
  \end{cases}
$$
We prove that $g$~is $TT_M$ which will be a contradiction.
So let $\tau$~be an $M$-tension on~$G$, we are to prove that $\tau g$
is an $M$-tension on~$F$. By the choice of~$l$, graph $f(F) \cap G'$
doesn't contain an $M$-unbalanced circuit (there is no that short
unbalanced circuit in~$G'$), hence any constant mapping is an $M$-tension.
So we may choose a tension $\tau '$ on~$G+G'$
that equals a constant~$\tau(e_0)$ on~$f(F)\cap G'$ and extends~$\tau$.
Clearly $\tau g$ is the same function as $\tau' f$, hence it is a tension.

For the last part of statement of the theorem, $G \hom G+G' \hom H$
follows immediately (using Lemma~\ref{SIL}, part 1). If we had
$H \hom K$ or $K \hom G$, then by Lemma~\ref{homo} the homomorphism
induces a $TT_M$ mapping $H \TT[M] K$ (or $K \TT[M] G$, respectively), 
a contradiction.
\qed
\end {proof}

To state Theorem~\ref{AntiExtInt} in a concise form
we define open and closed intervals in order $\lTT$. 
Let $(G,H)_M = \{ G' \mid   G \lTT_M G' \lTT_M H\}$
and $[G,H]_M= \{ G' \mid  G \leTT_M G' \leTT_M H\}$.
Similarly, define $(G,H)_h$ and $[G,H]_h$---intervals
in order $\lh$. Lemma~\ref{homo} implies that 
$[G,H]_h \subseteq [G,H]_M$ for any group~$M$. 
On the contrary, none of the two possible inclusions between
$(G,H)_h$ and~$(G,H)_M$ is valid for every $G$, $H$.
Therefore the additions in the following corollaries do indeed provide a
strengthening, we will use this strengthening in Section~\ref{sec:bddac}.

\begin {corollary}   \label{density}
Suppose $G \lTT_M H$ and $H \nTT[M] \edge$.
Then $(G,H)_M$ is nonempty. 
If in addition $G \lh H$ then $(G,H)_M \cap (G,H)_h$ is nonempty.
\end {corollary}

\begin {corollary}   \label{AntichainExt}
Suppose $G \lTT_M H$ and $H \nTT[M]\edge$.
Then any finite antichain of $\lTT_M$ restricted to~$(G,H)$
can be extended. 
If in addition $G \lh H$ then any finite antichain of $\lTT_M$
restricted to $(G,H)_M \cap (G,H)_h$ can be extended.
\end {corollary}

\begin{remark}
Throughout this section we need to assume $H \gTT_M \edge$:
for example in Corollary~\ref{density} there is no graph~$K$
satisfying $K_1 \lTT_M K \lTT_M \edge$ (if $K$ has no edge then
it maps to $K_1$, otherwise $\edge$ maps to it). We may say
that $(K_1,\edge)$ is a \emph{gap}.

If $M = \zet_2$ all results of this section hold for undirected
graphs, too, as all orientations of an undirected graph are
$TT_2$-equivalent.
\end{remark}

\begin{remark}
If $M$ is finite we can prove Lemma~\ref{SIL} easily 
by using the construction $\Delta(G)$ (and its variant for 
general group~$M$). For details, see~\cite{NS}.
\end{remark}

\section{Universality of $TT_2$ order}
\label{sec:universality}

In this section we restrict our attention to $TT_2$ mappings 
and consequently to undirected graphs. We first construct a
particular $TT_2$-rigid graph. (By Corollary~\ref{TTrigid} such graph
exists, but we need some additional properties.) Then we use this graph to
provide a faithful functor from the category of homomorphisms to the
category of $TT_2$~mappings.

\begin{figure}
\epsfbox{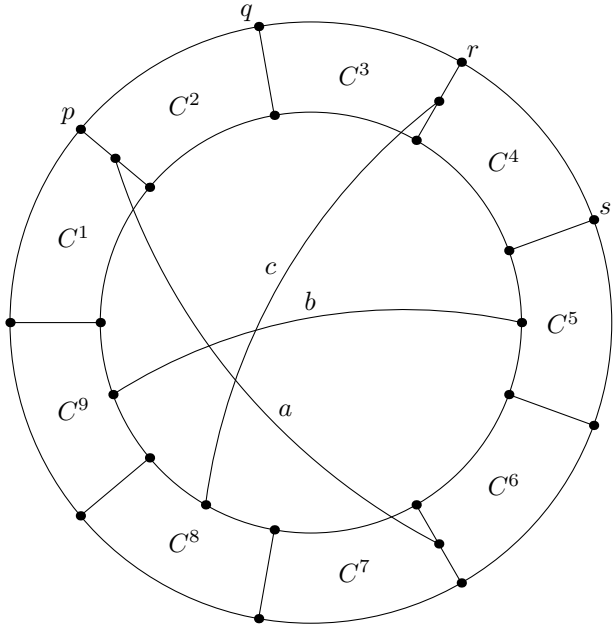}
\caption{A $TT_2$-rigid graph}
\label{fig:rigid}
\end{figure}

\begin {lemma}   \label{smallrigid}
Let~$S$ be the graph in~Figure~\ref{fig:rigid}.
\begin {enumerate}
\item $S$ is $TT_2$-rigid, i.e.\ the only $TT_2$~mapping $S \to S$
  is the identity.
\item Suppose $G$ is a graph that contains edge-disjoint copies
  of~$S$: $S_1$, \dots, $S_t$. Suppose $G$~does not contain
  triangles nor pentagons, except those pentagons that are contained 
  in some~$S_i$. Then the only $TT_2$~mapping $S \to G$ is the 
  identity mapping to some~$S_i$.
\end {enumerate}
\end {lemma}

\begin {proof}
We will prove the second part, which implies the first (by taking
$G=S$). Consider a $TT_2$~mapping $f:S \to G$.
Let pentagons in~$S$ be denoted $C^1$, \dots, $C^9$ as in the
figure, note that there are no other pentagons in $S$.
As there are no triangles in $G$ and the only pentagons are contained
in some $S_k$, we can deduce by Lemma~\ref{altdef} that each $C^i$
maps to a pentagon in some $S_k$ (possibly different $k$ for different~$i$).

Pentagon $C^i$ shares an edge with $C^j$ iff $i$~and $j$ differ
by~1 (modulo 9). As sharing an edge is preserved by any mapping
and since different copies of~$S$ in~$G$ are edge-disjoint, we conclude
that there is a copy of~$S$ in~$G$ (to simplify the notation, we
will identify this copy with~$S$) and a bijection
$p:[9] \to [9]$ such that $f(C^i)= C^{p(i)}$ for each~$i$;
moreover $p$~preserves the cyclic order.
Next we note that the size of the intersection of neighbouring pentagons
is preserved too. There are exactly three pairs of pentagons that
share two edges: $\{C^1,C^2\}$, $\{C^3,C^4\}$, $\{C^6,C^7\}$. As
the pairs $\{C^1,C^2\}$~and $\{C^3,C^4\}$ are adjacent,
the pairs $\{C^5,C^6\}$~and $\{C^3,C^4\}$ have a common neighbouring
pentagon, while the pairs $\{C^5,C^6\}$~and $\{C^1,C^2\}$ do not, we see
that $p$ is the identity; that is $f(C^i) = C^i$ for each~$i$.

We still have to prove that $f$~does not permute
edges in the respective pentagons. Let~$C^o$ be the outer cycle
and note it is the only 9-cycle in~$S$ that shares exactly one
edge with each~$C^i$. Hence, $f$ is an identity on~$E(C^o)$.
This means that $f$ can only permute two edges that share
an endpoint of some of the edges~$a$,~$b$, and~$c$.

Edge~$a$ is a part of a 7-cycle $C^a$ that has four edges in common
with~$C^o$. Now, $C^o$ is preserved by~$f$, and there is no other 7-cycle
in~$S$ with the same intersection with $C^o$. Thus, $C^a$ is preserved as
well, in particular $a$~and the edges incident with it are preserved.
Edge~$b$ is a part of a 7-cycle $C^b$ that intersects $C^5$,~$C^6$, $C^7$,
$C^8$ and~$C^9$. Since the edges it has in common with $C^6$, $C^7$, and~$C^8$ are
preserved by~$f$ (at least set-wise), and there is no other 7-cycle including 
these edges, $C^b$ is preserved too, in particular $b$~and the edges incident 
with it are preserved. Similarly, $c$ is contained in an $8$-cycle that has five of
its edges fixed, hence it is fixed by~$f$.
\qed
\end {proof}

\begin {theorem}   \label{faithful}
There is a mapping $F$ that assigns (undirected) graphs to
(undirected) graphs, such that for any graphs $G$,~$H$ (we stress that
we consider loopless graphs only) holds
$$
   G \hom H  \iff   F(G) \TT[2] F(H) \,.
$$
Moreover $F$ can be extended on mappings between graphs:
if $f:G \to H$ is a homomorphism, then $F(f):F(G)\to F(H)$
is a $TT$~mapping and any $TT$~mapping
between $F(G)$ and~$F(H)$ is equal to~$F(f)$ for some 
homomorphism~$f : G \hom H$.
(In category-theory terms $F$ is an embedding of the category of
all graphs and their homomorphisms into the category of
all graphs and all $TT_2$-mappings between them.)
\end {theorem}

\begin{figure}
\epsfbox{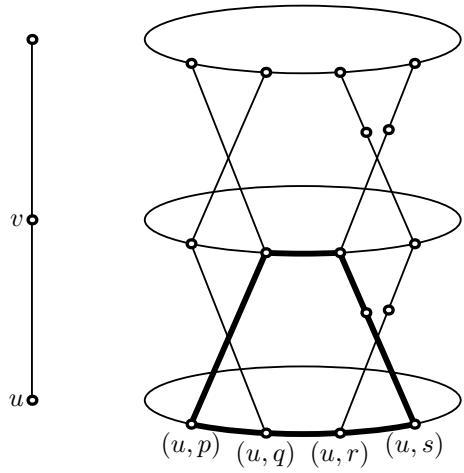}
\caption{Example of construction of $F(G)$ for $G = P_2$.
The 7-cycle used in the proof of Theorem~\ref{faithful} is drawn bold.}
\label{fig:faithful}
\end{figure}

\begin {proof}
Let $S$ be the graph from Lemma~\ref{smallrigid}, let 
$p$, $q$, $r$, $s$ be its vertices as denoted in
Figure~\ref{fig:rigid}.
For a graph~$G$, let the vertices of~$F(G)$ be
$(V(G)\times V(S)) \mathbin{\dot{\cup}} (E(G) \times \{1,2\})$. 
On each set $\{v\}\times V(S)$ we place a copy of~$S$,
it will be denoted by~$S_v$. For an edge~$uv$ of~$G$ we introduce
edges~$(u,p)(v,q)$, $(u,q)(v,p)$ (we refer to them as 
to \emph{add-on edges})
and paths of length two from~$(u,r)$ to~$(v,s)$ and from $(u,s)$
to~$(v,r)$ (we refer to these as to add-on paths, the middle vertices of
these paths are $(uv,1)$ and~$(uv,2)$).  There are no other edges
in~$F(G)$. See Figure~\ref{fig:faithful} for an example of the
construction.
As we wish to apply Lemma~\ref{smallrigid}, we first show that
$F(G)$~contains no triangles and only those pentagons that are 
contained in some $S_v$. Suppose $C$~is a cycle violating this.
If $C$ contains some add-on path, it is easy to check
that the length of $C$ is at least six. If it is not then $C$~has to
contain some add-on edges (as $S$ is triangle-free).  If it contains
only add-on edges and copies of the edge $pq$ then it has even length;
otherwise it has length at least seven.

It is clear how to define $F(f)$ for a homomorphism 
$f: G \to H$---$F(f)$ maps each~$S_v$ in~$G$ to~$S_{f(v)}$ in~$H$ in the
only way, the edges between different copies of~$S$ are mapped in the
``canonical'' way. Clearly $F(f)$ is a $TT$~mapping induced by a
homomorphism.

The only difficult part is to show, that for every $g: F(G) \TT F(H)$
there is an $f: G \hom H$ such that $g = F(f)$.
So let $g$ be such a mapping. By Lemma~\ref{smallrigid}
each copy of~$S$ is mapped to a copy of~$S$, to be precise, there
is a mapping $f: V(G) \to V(H)$ such that $g$ maps $S_v$ to
$S_{f(v)}$. Let $uv$ be an edge of~$G$.
First, we show that $f(u) \ne f(v)$. Suppose the contrary and
consider the 7-cycle $(u,p)$, $(u,q)$, $(u,r)$, $(u,s)$, $x$, 
$(v,r)$, $(v,q)$ ($x$ is the middle vertex of an add-on path).
Since $S$ is rigid, edges $(u,q)(u,r)$
and $(v,q)(v,r)$ map to the same edge, hence the algebraical image of the
other five edges is a cycle. However, there is no cycle of length at
most five containing edges $pq$ and $rs$, a contradiction.

Considering again the image of the same cycle shows that $f(u)$~and
$f(v)$ are connected by an edge of~$H$, which finishes the proof.
\qed
\end {proof}

\begin{remark}
It is interesting to note that graphs $F(G)$ are all triangle-free. We believe
that the construction from Theorem~\ref{faithful} can be modified to work 
for other groups than~$\zet_2$, some modification can possibly produce even 
graphs of girth at least~$g$, for any given~$g$.
If we consider graphs containing complete graphs then the situation becomes
easier. In fact (as we show in the next section), $TT$ mappings coincide
with homomorphisms on a large class of graphs (called nice graphs)---see
Theorem~\ref{niceTT} and the discussion below it.
\end{remark}

\section{Random graphs and $TT$~mappings}
\label{sec:random}

In this section we investigate cut-continuous, i.e.\ $TT_2$ mappings, only; 
that is we restrict our attention to the case $M=\zet_2$ and to
undirected graphs. We study whether typically (in the sense of random
graphs) a $TT_2$~mapping is induced by a homomorphism. Recall, that a
graph~$G$ is said to be \emph{homotens} if for any graph~$H$ any
$TT_2$ mapping $G \TT[2] H$ is induced by a homomorphism. The
main result of this section is that most graphs are homotens.

We consider the random graph model $\Ge_n$, that is every 
graph with vertices $\{1, 2, \dots, n\}$ has the same probability
(although some of the results can be modified for other models too).
As it is usual in the random graph setting, we study whether some
graph property $P$ holds almost surely (a.s.), that is whether 
$$
      \lim_{n \to \infty} Pr_{G \in \Ge_n}[\mbox{$G$ has $P$}] = 1 \,.
$$

We start with a useful notion that will help us to handle $TT$ mappings
(see Theorem~\ref{niceTT}).
We call a graph~$G$ \emph{nice} if the following holds
\begin {enumerate}
  \item every edge of~$G$ is contained in some triangle
  \item every triangle in~$G$ is contained in some copy of~$K_4$
  \item every copy of~$K_4$ in~$G$ is contained in some copy of~$K_5$
  \item for every $K$, $K'$ that are copies of~$K_4$ in~$G$ there is a 
    sequence of vertices
    $v_1$, $v_2$, \dots, $v_t$ such that
    \begin {itemize}
      \item $V(K) = \{v_1, v_2, v_3, v_4\}$,
      \item $V(K') = \{v_{t}, v_{t-1}, v_{t-2}, v_{t-3}\}$, 
      \item $v_i v_j$ is an edge of~$G$ whenever
        $1 \le i < j \le t$ and $j \le i+3$.
    \end {itemize}
\end {enumerate}

\begin {lemma}   \label{TTK5}
Let $f: K_5 \TT H$, where $H$ is any loopless graph.
Then $f$ is induced by an injective homomorphism (that is, 
by an embedding). Moreover, this isomorphism is uniquely determined.
\end {lemma}

\begin {proof}
Suppose $f(K_5)$ is a four-colorable graph. A composition 
of $TT$ mapping $f : K_5 \to f(K_5)$ with a $TT$ mapping
induced by a homomorphism $f(K_5) \to K_4$ gives
$K_5 \TT K_4$. Consider three cuts of size 4 in $K_4$; they
cover every edge exactly twice. Hence, their preimages are three cuts
in~$K_5$ that cover every edge exactly twice. But $K_5$ has 20 edges,
while the largest cut has only $2\cdot 3=6$~edges.

Hence, chromatic number of $f(K_5)$ is at least five. As it has
at most 10 edges, the chromatic number is exactly five. Let 
$V_1$, \dots, $V_5$ be the color classes. There is exactly
one edge between two distinct color classes (otherwise
the graph is four-colorable). Hence, $f$ is a bijection.
Next, $|V_i| = 1$ for every~$i$ 
(as otherwise, we can split one color-class to several 
pieces and join these to the other classes; again, the
graph would be four-colorable). Consequently, $f(K_5)$
is isomorphic to~$K_5$.

We call star a set of edges sharing a vertex. We know that preimage
of every star is a star, hence as $f$ is a bijection, also image of
every star is a star. Stars sharing an edge map to stars sharing an
edge, hence $f$ is induced by a homomorphism.
\qed
\end {proof}

\begin {theorem}   \label{niceTT}
Let $G$ be a nice graph, let $f : G \TT H$. Then $f$ is induced by a
homomorphism. Shortly, every nice graph is homotens.
\end {theorem}

\begin {proof}
Let $K$ be a copy of~$K_5$ in~$G$. 
By Lemma~\ref{TTK5} the restriction of~$f$ to~$K$
is induced by a homomorphism, let it be denoted by~$h_K$. 
If $K$ is a copy of~$K_4$ in~$G$, by the third condition
from the definition of nice it is contained in some $K'$---copy
of~$K_5$. The restriction $h_K = h_{K'} |_K$ induces~$f$ on~$K$;
clearly such $h_K$ is unique (it does not depend on the
choice of~$K'$). 

As every edge is contained in some copy of~$K_4$, 
it is enough to prove that there is a common extension
of all homomorphisms $\{h_K \mid K \subseteq G,\ K \simeq K_4\}$
(we may define it arbitrarily on the isolated vertices of~$G$).

We say that $h_K$ and~$h_{K'}$ \emph{agree} if for
any $v \in V(K) \cap V(K')$ we have $h_K(v) = h_{K'}(v)$. 
Thus, we need to show that any two homomorphisms $h_K$, $h_{K'}$
($K\simeq K'\simeq K_4$) agree.

Let first $K$, $K'$ be copies of $K_4$ that intersect in a triangle. 
Then $h_K$ and~$h_{K'}$ agree (note that this does not necessarily
hold if the intersection is just an edge). 

Now suppose $K$, $K'$ are copies of $K_4$ that have a common vertex~$v$. 
Since $G$ is nice, we find $v_1$, $v_2$, \dots, $v_t$ as in the
definition. Let $K_i = G[\{v_i, v_{i+1}, v_{i+2}, v_{i+3}\}]$; every 
$K_i$ is a copy of $K_4$, $K_1 = K$ and $K_{t-3}=K'$.
Suppose $v = v_l = v_r$, where $l \in \{1, 2, 3, 4\}$, 
$r \in \{t-3, t-2, t-1, t\}$. 
Consider a closed walk $W = v_{l}, v_{l+1}, \dots, v_{r-1}, v_r$.
Let $v'_i = h_{K_i}(v_i)$ for $l \le i \le r-3$
and $v'_i = h_{K_{r-3}}(v_i)$ for $r-3 \le i \le r$. 
Homomorphisms $h_{K_i}$ and $h_{K_{i+1}}$ agree, hence
$v'_i v'_{i+1} = f(v_i v_{i+1})$ is an edge of~$H$. So
$W' = v'_{l}, v'_{l+1}, \dots, v'_{r-1}, v'_r$.

Let $\phi(e)$ be the number of occurrences of~$e$ in~$W$ taken modulo 2. 
Clearly $\phi$ is a $\zet_2$-flow. 
Similarly, define $\phi'(e)$ as the number of occurrences of~$e$ in~$W'$
taken modulo 2. We have $\phi' = \phi_f$, hence by Lemma~\ref{altdef}
$\phi'$ is a flow. This can happen only if $W'$ is a closed walk, 
that is $v'_l = v'_r$. 

By definition, $v'_r = h_{K'}(v)$. As mappings $h_{K_i}$ and $h_{K_{i+1}}$
agree, we have that $h_{K_i}(v_{i+3}) = h_{K_{i+3}}(v_{i+3})$.
Consequently, $v'_l = h_K(v)$, which finishes the proof.
\qed
\end {proof}

Let us remark that Theorem~\ref{niceTT} may be used to prove
Theorem~\ref{faithful} in a different way. To do this, it suffices to
modify  the replacement operation (\cite{HN}) in such a way that the
resulting graph $F(G)$ is nice. (See \cite{HN} for a nice example of a
nice rigid graph.)

Consider the countable random graph $\Ge_{\omega}$. Surprisingly, 
it is almost surely isomorphic to a particular graph, the so-called
Rado graph. This is a remarkable graph (it is homogeneous and
it contains every countable graph as an induced subgraph), 
see~\cite{Cameron} for more detailed discussion.

\begin {lemma}   \label{nicegraphs}
Random graph from $\Ge_n$ is almost surely nice.
The Rado graph is nice. 
\end {lemma}

\begin {proof}
We prove the first statement, the second is proved in exactly
the same way, except we do not have to take the limit.

For $S \subseteq V(G)$ (where $G = \Ge_n$)
write~$C_S$ for the event, that there is a common neighbor
for all vertices in~$S$. If $|S|=4$, clearly the probability of~$C_S$
is~$(1-\frac 1{2^{s}})^{n-s}$. 
As $\binom ns \cdot (1-\frac 1{2^{s}})^{n-s}$ tends to zero for any
fixed $s$, $C_S$ holds a.s. for all $S$ with size at most 4.
This implies the first three conditions on~$G$.

To prove the last condition, let $K$, $K'$ be two copies of $K_4$.
Denote vertices of $K$ by $v_1$, $v_2$, $v_3$, $v_4$, and
vertices of $K'$ by $v_8$, $v_9$, $v_{10}$, $v_{11}$ (in any order).
If we find a triangle that is connected to every vertex in $K \cup K'$, 
we may denote its vertices by $v_5$, $v_6$, $v_7$ and we are done.
For a given three-element set $S \subseteq V(G) \setminus (V(K)\cup V(K'))$
the probability that $S$ induces a triangle and is connected
to all vertices in $V(K)\cup V(K'))$ is at least $2^{-21}$, hence
the probability that there is no such $S$ is at most
$(1-2^{-21})^{n-6/3}$. As the number of possible pairs $(K,K')$
is~$O(n^8)$, this concludes the proof.
\qed
\end {proof}

\begin {lemma}   \label{nicecomplete}
The complete graph~$K_n$ is nice whenever $n \ge 5$.
\end {lemma}

\begin {proof}
The straightforward verification is left to the reader.
\qed
\end {proof}

From Theorem~\ref{niceTT}, Lemma~\ref{nicegraphs},
and Lemma~\ref{nicecomplete} we immediately get the
following corollary (a different proof of
$K_4 \lTT K_5 \lTT \cdots$ is given in~\cite{LMT}).

\begin {corollary}   \label{homotens}
\begin {enumerate}
  \item Random graph from~$\Ge_n$ is almost surely homotens.
  \item The Rado graph is homotens.
  \item The complete graph~$K_n$ is homotens whenever $n \ge 5$.
    In particular, in the $TT_2$ order we have
    $$
      K_3 \eqTT K_4 \lTT K_5 \lTT K_6 \lTT K_7 \lTT \cdots \,.
    $$
\end {enumerate}
\end {corollary}

Corollary~\ref{homotens} enables us to prove a $TT$ version
of the following result about homomorphisms of random graphs. 
(The original theorem appears in~\cite{KR}, see also Section~3.6
of~\cite{HN}.)

\begin {theorem}[\cite{KR}]   \label{homorigid}
Random graph is almost surely rigid (with respect to homomorphism).
There are 
$$
  \frac{1}{n!} \binom { \binom n2 }{ \floor{\frac 12 \binom n2}} (1-o(1))
$$
graphs on~$n$~vertices with no homomorphism between any two of them
and with only identical homomorphism on each of them.
\end {theorem}

\begin {corollary}   \label{TTrigid}
  Random graph is almost surely $TT$-rigid. There are 
$$
  \frac{1}{n!} \binom { \binom n2 }{ \floor{\frac 12 \binom n2}} (1-o(1))
$$
  pairwise $TT$-incomparable $TT$-rigid graphs on $n$~vertices. 
\end {corollary}


\begin{remark}
The method of this section may be used for~$M\ne \zet_2$ as well. 
In fact, if $M$ is not a power of~$\zet_2$, 
we can prove analogy of Lemma~\ref{TTK5} for $K_4$. 
Then we can prove stronger version of Theorem~\ref{niceTT}---for
any group~$M$, a nice graph is $M$-homotens; we can
even slightly weaken the definition of ``nice'' if $M$ is not a power
of~$\zet_2$. Similarly, we can generalize other results of this
section. For details, see~\cite{NS} and~\cite{RSthesis}.
\end{remark}

\section{Influence of the group}
\label{sec:group}

In this section we study how the notion of $M$-tension-continuous 
mapping depends on the group~$M$. 
Although the existence of $M$-tension-continuous mappings seems
to be strongly dependent on the choice of~$M$ we prove here 
(in Theorem~\ref{numberofTT}) that this dependence
relates only to the cyclical structure of~$M$.

Throughout this section, $G$, $H$ will be two graphs, 
$f: E(G) \to E(H)$ a mapping, and $M$, $N$ groups, recall we consider only
abelian groups (as is usual in the study of group-valued flows).
As we are interested mainly in finite graphs, we can restrict our
attention to finitely generated groups---clearly $f$ is
$M$-tension-continuous iff it is $N$-tension-continuous for
every finitely generated subgroup of~$M$.

Hence, we can use the classical characterization of finitely
generated Abelian groups (see e.g.~\cite{Lang}).
\begin {theorem}   \label{abgroups}
For a finitely generated abelian group $M$ there are integers
$\alpha$, $k$, $\beta_i$, $n_i$ ($i = 1, \dots, k$)
so that
\begin{equation}\label{eq:M}
  M \simeq 
     \zet^\alpha \times \prod_{i=1}^k \zet_{n_i}^{\beta_i} \,.
\end{equation}
\end {theorem}

For a group $M$ in the form~(\ref{eq:M}), denote $n(M) = \infty$ if
$\alpha > 0$, otherwise let $n(M)$ be the least common multiple of 
$\{n_1, \dots, n_k\}$. 

As a first step to complete characterization we consider a  
specialized question: given a $TT_M$~mapping, when can we conclude
that it is $TT_N$ as well?

\begin {lemma}   \label{MN}
\begin {enumerate}
  \item If $f$ is $TT_\zet$ then it is $TT_M$ for any $M$.
  \item Let $M$ be a subgroup of $N$. If $f$ is $TT_N$ then it
    is $TT_M$.
\end {enumerate}
\end {lemma}

\begin {proof}
1. This appears as Theorem~4.4 in~\cite{DNR}.

2. Let $\tau$ be an $M$-tension on~$H$. As $M \le N$, we may
regard $\tau$ as an $N$-tension, hence $\tau f$ is an $N$-tension
on~$G$. As it attains only values in the range of $\tau$, hence
in~$M$, it is an $M$-tension, too.
\qed
\end {proof}

\begin {lemma}   \label{M12}
Let $M_1$, $M_2$ be two abelian groups.
Mapping $f$ is $TT_{M_1}$ and $TT_{M_2}$ if and only if
it is $TT_{M_1 \times M_2}$.
\end {lemma}

\begin {proof}
As $M_1$, $M_2$ are subgroups of $M_1 \times M_2$, one implication
follows from the second part of Lemma~\ref{MN}. For the other implication
let $\tau$ be an $M_1 \times M_2$ tension on~$H$. Write 
$\tau = (\tau_1, \tau_2)$, where $\tau_i$ is an $M_i$-tension on~$H$.
By assumption, $\tau_i f$ is an $M_i$ tension on~$G$, consequently
$\tau f = (\tau_1 f, \tau_2 f)$ is a tension too.
\qed
\end {proof}

The following (somewhat surprising) lemma shows that we can restrict
our attention to cyclic groups only.

\begin {lemma}   \label{onlycyclic}
\begin {enumerate}
  \item If $n(M) = \infty$ then $f$ is $TT_M$ if and only if it is $TT_\zet$.
  \item Otherwise $f$ is $TT_M$ if and only if it is $TT_{n(M)}$.
\end {enumerate}
\end {lemma}

\begin {proof}
By previous lemmas. Note that $\zet_n$ is a subgroup 
of~$\prod_{i=1}^k \zet_{n_i}^{\beta_i}$.
\qed
\end {proof}

By a theorem of Tutte (see~\cite{Diestel}), the number of nowhere-zero flows on 
a given graph does depend only on the size of the group (that is, surprisingly, 
it does not depend on the structure of the group). 
Before proceeding in the main direction of this section, let us
note a consequence of Lemma~\ref{onlycyclic}, which is an analogy of the 
Tutte's theorem.

\begin {theorem}   \label{numberofTT}
Given graphs~$G$, $H$, the number of $TT_M$ mappings from~$G$ to~$H$ 
depends only on~$n(M)$.
\end {theorem}

Lemma~\ref{onlycyclic} suggests to define for two graphs the set
$$
  TT(G,H) = \{ n \ge 1 \mid \hbox{there is $f: E(G) \to E(H)$ 
          such that $f$ is $TT_n$} \}
$$
and for a particular $f: E(G) \to E(H)$ 
$$
  TT(f,G,H) = \{ n \ge 1 \mid \hbox{$f$ is $TT_n$} \}  \,.
$$

Remark that most of these sets contain 1: $\zet_1$ is a trivial 
group, hence any mapping is $TT_1$. 
Therefore $1 \in TT(f,G,H)$ for every $f: E(G) \to E(H)$, 
while $1 \in TT(G,H)$ iff there exists a mapping $E(G) \to E(H)$.
This happens precisely when $E(H)$ is nonempty or $E(G)$ is empty.

\begin {lemma}   \label{TTf-fin}
Either $TT(f,G,H)$ is finite or $TT(f,G,H)=\en$.
In the latter case $f$ is $TT_\zet$.
\end {lemma}

\begin {proof}
It is enough to prove that $f$ is $TT_\zet$ if it is $TT_n$
for infinitely many integers~$n$. To this end, take a $\zet$-tension 
$\tau$ on~$H$. As $\tau_n : e \mapsto \tau(e) \bmod n$ is a
$\zet_n$-tension, $\tau_n f = \tau f \bmod n$ is a $\zet_n$-tension
whenever $f$~is $TT_n$.
To show $\tau$ is a $\zet$-tension consider a circuit~$C$ and let
$s$ be the ``$\pm$-sum'' (in~$\zet$) along~$C$. As $s \bmod n = 0$
for infinitely many values of~$n$, we have $s=0$.
\qed
\end {proof}

Any $f$ induced by a homomorphism provides an example where
$TT(f,G,H)$ is the whole $\en$. For finite sets, the situation is more
interesting. By the next theorem the sets $TT(f,G,H)$
are precisely ideals in the divisibility lattice.

\begin {theorem}   \label{TTf-str}
Let $T$ be a finite subset of~$\en$. Then the following are
equivalent.
\begin {enumerate}
  \item There are $G$, $H$, $f$ such that $T = TT(f, G, H)$. 
  \item There is $n \in \en$ such that $T$ is the set of 
    all divisors of~$n$.
\end {enumerate}
\end {theorem}

\begin {proof}
First we show 1. implies 2.
The set $T$ has the following properties
\begin {enumerate}[(i)]
  \item If $a \in T$ and $b | a$ then $b \in T$. 
    (We use the second part of Lemma~\ref{MN}: if $b$ divides $a$, then
    $\zet_b \le \zet_a$.)
  \item If $a, b \in T$ then the least common multiple of $a$, $b$
     is in~$T$. 
    (We use Lemma~\ref{MN} and Lemma~\ref{M12}: if 
    $l = \lcm(a,b)$ then $\zet_l \le \zet_a \times \zet_b$.)
\end {enumerate}

Denote $n$ the maximum of~$T$. By (i), all divisors of~$n$ are in~$T$.
If there is a $k \in T$ that does not divide~$n$ then $\lcm(k,n)$
is element of~$T$ larger than~$n$, a contradiction.

For the other implication, let $f$ be the only mapping from~$\orC_n$
to~$\edge$. Then $TT(f,\orC_n,\edge) = T$: mapping~$f$
is~$TT_k$ iff for any $a \in \zet_k$ the constant mapping
$E(\orC_n) \mapsto a$ is a $\zet_k$-tension; this occurs
precisely when $k$~divides~$n$. 
\qed
\end {proof}

Let us turn to description of sets $TT(G,H)$.
Although we are working with finite graphs throughout the paper, we
stress here that $G$, $H$ are finite graphs---in contrary with most
of other results, this one is not true for infinite graphs.

\begin {lemma}   \label{TT-fin}
Let $G$, $H$ be \emph{finite} graphs.
Either $TT(G,H)$ is finite or $TT(G,H)=\en$.
In the latter case $G \leTT_\zet H$.
\end {lemma}

\begin {proof}
As in the proof of~Lemma~\ref{TTf-fin}, the only difficult
step is to show that if $G \leTT_n H$ for infinitely many 
values of~$n$, then $G \leTT_\zet H$. As $G$~and~$H$ are finite, 
there is only a finite number of possible mappings between
their edge sets. Hence, there is one of them, say $f: E(G) \to E(H)$, 
that is $TT_n$ for infinitely many values of~$n$. By
Lemma~\ref{TTf-fin} we have $f: G \TT[\zet] H$.
\qed
\end {proof}

When characterizing the sets $TT(G, H)$ we first remark that
an analogue of Lemma~\ref{M12} does not hold: there is a 
$TT_M$~mapping from~$\orC_9$ to~$\orC_7$ for
$M=\zet_2$ (mapping induced by a homomorphism of the 
undirected circuits) and for $M=\zet_3$ (e.g. a constant mapping), 
but not the same mapping for both, hence there
is no $TT_{\zet_2\times\zet_3}$~mapping. We will see that
the sets $TT(G,H)$ are precisely down-sets in the divisibility poset.
First, we prove a lemma that will help us to construct pairs of graphs
$G$, $H$ with a given $TT(G,H)$.
Integer cone of a set $\{s_1, \dots, s_t\} \subseteq \en$ is
the set $\{ \sum_{i=1}^t a_i s_i \mid a_i \in \zet, a_i \ge 0 \}$.

\begin {lemma}   \label{circuits}
Let $A$, $B$ be non-empty subsets of~$\en$, $a \in \en$, define
$G = \bigcup_{a \in A} \orC_a$, and~$H = \bigcup_{b \in B} \orC_b$.
Then there is a $TT_n$~mapping from~$G$ to~$H$
if and only if 
$$
  \hbox {$A$ is a subset of the integer cone of $B \cup \{ n\}$} \,.
$$
\end {lemma}

\begin {proof}
We use Lemma~\ref{altdef}.
Consider a flow $\phi_a$ attaining value 1 on $\orC_a$ and 0
elsewhere. Algebraical image of this flow is a flow, hence
it is (modulo $n$) a sum of several flows along the cycles
$\orC_b$, implying $a$~is in integer cone of~$B \cup \{n\}$.
On the other hand if $a = \sum_i b_i + cn$ then we can
map any $c$ edges of~$\orC_a$ to one (arbitrary) edge of~$H$, 
and for each $i$ any (``unused'') $b_i$ edges bijectively to
$\orC_{b_i}$. After we have done this for each $a \in A$ we will have
constructed a $\zet_n$-tension-continuous mapping from~$G$ to~$H$.
\qed
\end {proof}

\begin {theorem}   \label{TT-str}
Let $T$ be a finite subset of~$\en$. Then the following are
equivalent.
\begin {enumerate}
  \item There are $G$, $H$  such that $T = TT(G, H)$. 
  \item There is a finite set $M \subset \en$ such that 
  $$
      T = \{ k \in \en ; (\exists m \in M) k | m \} \,.
  $$
\end {enumerate}
\end {theorem}

\begin {proof}
If $T$ is empty, we take $M$ empty. In the other direction, 
if $M$ is empty we just consider graphs such that $E(H)$ is empty and
$E(G)$ is not. Next, we suppose $M$ is nonempty.

By the same reasoning as in the proof of Theorem~\ref{TTf-str}
we see that if $a \in T$ and $b | a$ then $b \in T$.
Hence, 1 implies 2, as we can take $M = T$ (or, to make
$M$ smaller, let $M$ consist of the maximal elements of~$T$ in the
divisibility relation). 

For the other implication let $p > 4 \max M$ be a prime, 
let $p' \in (1.25p, 1.5p)$ be an integer. Let $A = \{p,p'\}$ and
$$
  B = \{ p-m; m \in M\} \cup \{ p'-m; m \in M\} \,;
$$
note that every element of~$B$ is larger than $\tfrac34 p$.
As in Lemma~\ref{circuits} we define
$G = \bigcup_{a \in A} \orC_a$, $H = \bigcup_{b \in B} \orC_b$.
We claim that $TT(G,H) = T$. By Lemma~\ref{circuits} it is immediate
that $TT(G,H) \supseteq T$. For the other direction take
$n \in TT(G,H)$. By Lemma~\ref{circuits} again, we can express
$p$ and~$p'$ in form
\begin{equation}\label{eq:formofp}
    \sum_{i=1}^t b_i + cn
\end{equation}
for integers $c, t \ge 0$, and $b_i \in B$.
\begin {itemize}
  \item If $t\ge 2$ then the sum in~(\ref{eq:formofp}) is at least $1.5p$;
     hence neither $p$ nor $p'$ can be expressed with $t\ge 2$.
  \item If $t=1$ then we distinguish two cases.
    \begin {itemize}[$\bullet$]
      \item $p = (p-m) + cn$, hence $n$ divides $m$ and $n\in T$.
      \item $p = (p'-m) + cn$, hence $p'-p \le m$. But
        $p'-p > 0.25p > m$, a contradiction.
    \end {itemize}
    Considering $p'$ we find that either $n \in T$ or
    $p' = (p-m) + cn$.
 \item Finally, consider $t=0$. If $p=cn$ then either $n=1 \in t$
   or $n=p$. (We don't claim anything about $p'$.)
\end {itemize}
To summarize, if $n \in TT(G,H) \setminus T$ then necessarily $n=p$.
For $p'$ we have only two possible expressions, $p'=cn$ and
$p' = (p-m) + cn$. We easily check that both of them lead to 
a contradiction. The first one contradicts $1.25 p < p' < 1.5 p$. 
In the second expression $c=0$ implies $p'<p$
while $c\ge 1$ implies $p'\ge 2p-m \ge 1.75p$, again a contradiction.
\qed
\end {proof}

\begin{remark}
This paper concentrates on $TT$~mappings. We remark, however, that 
the same proof yields a characterization of sets $XY(f,G,H)$ and $XY(G,H)$
for $XY \in \{FF, FT, TF\}$ (which are defined for $FF$, $FT$, and
$TF$ mappings in the same way as sets $TT(f,G,H)$ and
$TT(G,H)$ for $TT$ mappings).
\end{remark}

\section{Miscellanea}
\label{sec:misc}

\subsection{Complexity}
\label{sec:complexity}

\def\HOM{\mathop{\textsc{HOM}}}
\def\T{\mathop{\textsc{TT}}}

Let $\T_M(H)$ denote the problem of decision, whether for a given graph~$G$
there is a $TT_M$~mapping $G \TT[M] H$. The complexity of the related
problem $\HOM(H)$ (that is the testing of the existence of a homomorphism
to~$H$) is now well understood, at least for undirected graphs: $\HOM(H)$
is NP-complete if and only if $H$ contains an odd circuit, otherwise it
is in P (as it is equivalent to decide whether $G$ is bipartite), 
see~\cite{HN-complexity}.
In the same spirit, we wish to determine the complexity of the problem
$\T_M(H)$. 

\begin {theorem}   \label{cccompl}
Let $H$ be an undirected graph. Then $\T_{\zet_2}(H)$
is NP-complete if $H$ contains an odd circuit; otherwise
it is polynomial. 
\end {theorem}

\begin {proof}
By Theorem~\ref{cube}, problems $\T_2(H)$ and $\HOM(\Delta(H))$ have
the same answer for any graph~$G$, hence they have the same
complexity. Observe that $\Delta(H)$ is bipartite iff
$H$ is bipartite: $H$ and $\Delta(H)$ are $TT_2$ equivalent
and any graph is bipartite iff it admits a $TT_2$ mapping
to~$\edge$. Consequently, $\T_2(H)$ is NP-complete iff 
$H$~contains an odd circuit.
\end {proof}

For $M \ne \zet_2$ (or $\zet_2^k$), we may still reduce
$\T_M(H)$ to $\HOM(H')$ for a suitable graph~$H'$. However, now we deal 
with directed graphs, where the complexity of $\HOM$ is not characterized.
Another obstacle is that for $M=\zet$ the graph $H'$ is infinite.
(For $H$ infinite, the complexity of~$\HOM(H)$ was investigated
in~\cite{BN}.)

\subsection{Codes and $\chi/\chi_{TT}$}
\label{sec:codes}

In this section we first restate parts of~\cite{LMT} in our terminology.
Inspired by the definition of $\chi(G)$ via homomorphisms we may define 
$$
   \chi_{TT}(G) = \min \{ n ; G \TT[2] K_n \} \,.
$$
For random graphs, Corollary~\ref{homotens} implies that 
$\chi_{TT}(G) = \chi(G)$ almost surely. For general graph~$G$,
Lemma~\ref{homo} implies $\chi_{TT}(G) \le \chi(G)$, on the
other hand $\chi_{TT}(G) > \chi(G)/2$ follows from the
fact that homomorphisms and $TT_2$~mappings
to $K_{2^k}$ coincide (\cite{LMT}, \cite{DNR}).
More precise information on behaviour of $\chi(G)/\chi_{TT}(G)$ is
desirable. 

Consequently, let $\calG_n = \{G \mid G \TT[2] K_n \}$ and study
$\chi(G)$ for $G \in \calG_n$. By Lemma~\ref{cube}, $G \in \calG_n$
is equivalent to $G \hom \Delta(K_n)$. In other words, 
\begin {itemize}
  \item $\Delta(K_n) \in \calG_n$; and
  \item for every $G \in \calG_n$ we have $G \hom \Delta(K_n)$. 
\end {itemize}
This reduces the problem of behaviour of $\chi_{TT}(G)/\chi(G)$
to special values of~$G$.

\begin {problem}   \label{limitchi}
Determine the limit of $\chi(\Delta(K_n))/n$ (and in particular
decide, whether the limit exist). (We only know the fraction is always
in the interval $[1,2]$.)
\end {problem}

The chromatic number of $\Delta(K_n)$ was studied before
(with the same motivation) in~\cite{LMT}. 
In~\cite{HKSS}, the connection with injective chromatic number of hypercubes 
is presented. 
In~\cite{Dvorak} graphs $\Delta(K_n)$ are studied (as a special type
of graphs arising from hypercubes) in the context of embedding of
trees. It is claimed there that $\chi(\Delta(K_9)) \ge 13$.
There is a chapter on the topic in~\cite{JT} (``chromatic number of
cube-like graphs'').

If we see the vertices of~$\Delta(K_n)$
as $\{0,1\}^n$ then an independent set forms a ``code'' ---
a set where no two elements have Hamming distance 2. With some
more work we can use results from theory of error-correcting codes.
This approach was taken in~\cite{LMT} and~\cite{HKSS}. After 
using \cite{Best} they obtained what seems to be the
strongest result so far:
$\chi(\Delta(K_{n-3})) = n$ for $n = 2^k$ ($k\ge 2$).

We add a new piece of information to the picture:
if we restrict our attention to sparse graphs we see the same set of
values $\chi_{TT}(G)/\chi(G)$.

\begin {lemma}   \label{sparsechiTT}
Let $n$, $c$ be integers, $n \ge 3$. There is $G \in \calG_n$
such that $\chi(G) = \chi(\Delta(K_n))$ and $g(G) > c$. 
\end {lemma}

In the proof we will use Sparse incomparability lemma for homomorphisms in
the following form.

\begin {lemma}   \label{homoSIL}
Let $H$, $G_1$, \dots, $G_k$ be (undirected) graphs such that
$H$ is not bipartite and $H \nhom G_i$ for every $i$. Let $c$ be an
integer. Then there is an undirected graph~$G$ such that 
\begin {itemize}
  \item $g(G)>c$ (that is $G$~contains no circuit of size at most~$c$),
  \item $G \lh H$, and
  \item $G \nhom G_i$ for every~$i$.
\end {itemize}
\end {lemma}

\begin {proof}(of Lemma~\ref{sparsechiTT})
Suppose $\chi(\Delta(K_n)) = t$, hence $\Delta(K_n) \nhom K_{t-1}$. 
By Lemma~\ref{homoSIL} we get $G$ with $g(G) > c$ such that
$G \hom \Delta(K_n)$ and $G \nhom K_{t-1}$. Hence $G \in \calG_n$
and $\chi(G) > t-1$. On the other hand 
$\chi(G) \le \chi(\Delta(K_n)) = t$.
\qed
\end {proof}

\begin{remark}
In~\cite{LMT} it is proved that if we define $\chi_{TT_\zet}$
by $TT_\zet$ mappings, then $\chi_{TT_\zet}(G) = \chi(G)$ for
every graph~$G$. It may be worth to study $\chi_{TT_M}$ for
other groups~$M$, too.
\end{remark}

\subsection{Dualities in the $TT$ order}

\def\orT{\overrightarrow{T}\!}
\def\orP{\overrightarrow{P}\!}

Dualities were introduced as an example of good characterization which
can help to solve $\HOM(H)$ for some graphs~$H$.  We say that a tuple
$(F_1, \dots, F_t; H)$ forms a duality if for every $G$
$$
   G \hom H \iff (\forall i \in \{1, \dots, t\}) F_i \nhom G \,.
$$
It is well-known that $G$ has a homomorphism to~$\orT_n$
(transitive tournament with $n$~vertices)
iff it does not contain $\orP_{n+1}$ (path with $n+1$~vertices).
Hence, the pair $(\orP_{n+1};\orT_n)$ is a duality. 
If $(F_1, \dots, F_t; H)$ is a duality, we can solve
$\HOM(H)$ in polynomial time. Dualities are studied in a
sequence of papers, see~\cite{NT},~\cite{HN} and references there.
We present a sample of results:
\begin {itemize}
  \item for undirected graphs there are only trivial dualities
    $(K_2;K_1)$ and $(K_1;K_0)$.  
  \item for directed graphs, for any $t$ and any trees $F_1$, \dots, $F_t$, 
    there is an $H$ such that $(F_1, \dots, F_t; H)$ is a duality; there
    are no other dualities.
  \item similarly as for directed graphs, it is possible to characterize 
    all dualities for arbitrary relational systems. 
\end {itemize}

Here we adopt proof of the homomorphic case (for undirected graphs)
to characterize dualities
for $TT_M$, that is we characterize all tuples $(F_1, \dots, F_t; H)$ for which
\begin{equation}\label{eq:duality}
   G \TT H \iff (\forall i \in \{1, \dots, t\}) F_i \nTT G \,. 
\end{equation}
We suppose $M \ne \zet_1$ to avoid trivialities.

\begin {theorem}   \label{TTdualities}
For every group~$M$, there are no dualities in the $TT_M$ order,
up to the trivial ones, that is $H \eqTT_M K_1$ and for some~$i$ we
have $F_i \eqTT_M \edge$.
\end {theorem}

\begin {proof}
Let $(F_1, \dots, F_t; H)$ be a duality. Denote
$g = \max \{ g_{\sss M}(F_1), \dots, g_{\sss M}(F_t)\}$. If $g=\infty$, 
then there is an $i$ such that $F_i \TT[M] \edge$. In this case, 
the right-hand side of $(\ref{eq:duality})$ holds iff $G$~is edgeless. This 
is equivalent to $G \TT[M] H$ exactly when $H$~is edgeless, 
that is $H \eqTT_M K_1$.

If $g$ is finite, we consider a graph~$G$ such that $\chi(G) > c$
($c$ will be specified later) and all circuits in~$G$ are longer
than~$g$. (Such graphs exist by the celebrated theorem of Erd\H{o}s.)
We orient the edges of~$G$ arbitrarily. Now $F_i \nTT[M] G$ by
Lemma~\ref{invariant}, it remains to prove $G \nTT[M] H$.
So suppose the contrary;
by Lemma~\ref{onlycyclic} and~\ref{MN} we may suppose $M$~is finite.
By Theorem~\ref{cube} (and the remarks following it), there is a finite 
directed graph $H'$ such that $G \TT[M] H$ iff $G \hom H'$. 
Hence it is enough to choose $c = \chi(H')$.
\qed
\end {proof}

\subsection{Bounded antichains in the $TT$ order}
\label{sec:bddac}

In~\cite{DNR}, the following question is posed
(for $M=\zet_2$ as Problem~6.9, for $M = \zet$ implicitly at 
the end of Section~8).

\begin{problem*}[\cite{DNR}]
Is there an infinite antichain in the order $\leTT_M$, that
consists of graphs with bounded chromatic number?
\end{problem*}

Our approach provides a straightforward answer in a very 
strong form.

\begin {corollary}   \label{AntiChain}
For every $M$, there is an infinite antichain in the order
$\leTT_M$, that consists of graphs with chromatic number at most 3.
\end {corollary}

\begin {proof}

Let $G = \edge$ and choose a 3-colorable~$H$ such that $H \gTT_M \edge$:
we can take $H = \orC_3$ whenever $M$ is not a power of $\zet_3$.
In that case we choose $H = \orC_5$. 

Denote $I = (G,H)_h \cap (G,H)_M$. By Corollary~\ref{density}, $I$
is nonempty, hence choose $G_0 \in I$.
Now we inductively find (using Corollary~\ref{AntichainExt})
graphs $G_1$, $G_2$, \dots\ from~$I$ such that for every $k$, 
$G_0$, \dots, $G_k$ is an antichain in the order $\lTT_M$.
Hence $\{G_n, n\ge 0\}$ is an infinite
antichain, and as for every $i$, $G_i \hom H$, 
every $G_i$ is 3-colorable. 

For $M = \zet_2$, an alternative proof is provided 
by Theorem~\ref{faithful}: the homomorphism order is known to 
have infinite antichain of bounded chromatic number, this is mapped
to an infinite antichain in $\leTT_2$ of bounded $\chi_{TT_2}$, 
hence of bounded chromatic number ($\chi_{TT_2} \ge 2 \chi$, 
see Section~\ref{sec:codes}).
\qed
\end {proof}

\begin {remark}   
In the presented proof we can choose $H$ more carefully, 
namely we can let $H = \orC_p$ where $p$ is a large enough prime
(so that $\orC_p \nTT[M] \edge$). In this way, we obtain 
an infinite antichain of order $\lTT_M$ that consists of graphs
with circular chromatic number bounded by $2 + 1/\varepsilon$.
\end {remark}

\subsection{Differences between $\leTT_M$ and $\leh$}

We restrict our attention to $M = \zet_2$ and $M = \zet$, 
which seem to be the two most important cases.
As shown by Theorem~\ref{largeexamples} there are pairs
of arbitrary highly connected graphs $G$, $H$ such that
$G \eqTT_2 H$ while $G \lh H$. Note that this means that
$H$~is not nice: indeed, $H = \Delta(K_n)$ for a suitable~$n$, 
and although $\Delta(K_n)$ contains $K_n$, not every copy of~$K_4$ 
is contained in a copy of~$K_5$. 
On the contrary, Corollary~\ref{homotens} shows that for almost all graphs
$\leTT_2$ and $\leh$ coincide. It would be interesting to know, 
whether $\leTT_2$ and $\leh$ coincide for random regular graphs, 
or for sparse random graphs.

For $TT_\zet$ the situation is rather different. Any 
two oriented trees are $TT_\zet$-equivalent, hence we have
plenty of 1-connected graphs for which $\leTT_\zet$ and $\leh$ 
differ. For 2-connected examples, consider any permutation $\pi:
E(\orC_n) \to E(\orC_n)$. This is $TT_\zet$, but (except for $n$ of
them) is not induced by a homomorphism.  We may now use the
replacement operation of~\cite{HN}, that is we replace every edge
of~$\orC_n$ by a suitable graph (for every edge we use a different
graph). In this way we produce from the oriented circuit two graphs $G$
and~$H$, such that there is only one mapping $G \TT[\zet] H$, it
``obeys'' one of the permutations $\pi: E(\orC_n) \to E(\orC_n)$. So
if we choose $\pi$ that is not a cyclic shift, we obtain graphs such
that $G \eqTT_\zet H$ and $G \not \eqh H$. These graphs are (vertex)
2-connected, while they may have arbitrary edge-connectivity.
Presently, we do not know whether there are (vertex) 3-connected
graphs, where $\leTT_\zet$ and~$\leh$ differ; in fact we are not 
aware of (vertex) 3-connected graph that is not $\zet$-homotens.



\begin{thebibliography}{10}

\bibitem{Best}
Mark~R. Best and Andries~E. Brouwer, \emph{The triply shortened binary
  {H}amming code is optimal}, Discrete Math. \textbf{17} (1977), no.~3,
  235--245.

\bibitem{BN}
Manuel Bodirsky and Jaroslav Ne{\v{s}}et{\v{r}}il, \emph{Constraint
  satisfaction with countable homogeneous templates}, Computer science logic,
  Lecture Notes in Comput. Sci., vol. 2803, Springer, Berlin, 2003, pp.~44--57.

\bibitem{Cameron}
Peter~J. Cameron, \emph{The random graph}, The mathematics of Paul Erd{\H o}s,
  II, Algorithms Combin., vol.~14, Springer, Berlin, 1997, pp.~333--351.

\bibitem{DNR}
Matt DeVos, Jaroslav Ne{\v s}et{\v r}il, and Andr{\'e} Raspaud, \emph{On flow
  and tension-continuous maps}, KAM-DIMATIA Series \textbf{567} (2002).

\bibitem{Diestel}
Reinhard Diestel, \emph{Graph theory}, Graduate Texts in Mathematics, vol. 173,
  Springer-Verlag, New York, 2000.

\bibitem{Dvorak}
Tom{{\'a}}{{\v s}} Dvo{{\v r}}{{\'a}}k, Ivan Havel, Jean-Marie Laborde, and
  Petr Liebl, \emph{Generalized hypercubes and graph embedding with dilation},
  Proceedings of the 7th Fischland Colloquium, II (Wustrow, 1988), no.~39,
  1990, pp.~13--20.

\bibitem{Gue}
Bertrand Guenin, \emph{Packing {T}-joins and edge colouring in planar graphs},
  (to appear).

\bibitem{HKSS}
Ge{{\v n}}a Hahn, Jan Kratochv{{\'i}}l, Jozef {{\v S}}ir{{\'a}}{{\v n}}, and
  Dominique Sotteau, \emph{On the injective chromatic number of graphs},
  Discrete Math. \textbf{256} (2002), no.~1-2, 179--192.

\bibitem{HN-complexity}
Pavol Hell and Jaroslav Ne{\v{s}}et{\v{r}}il, \emph{On the complexity of
  {$H$}-coloring}, J. Combin. Theory Ser. B \textbf{48} (1990), no.~1, 92--110.

\bibitem{HN}
Pavol Hell and Jaroslav Ne{\v s}et{\v r}il, \emph{Graphs and homomorphisms},
  Oxford Lecture Series in Mathematics and Its Applications, Oxford University
  Press, 2004.

\bibitem{Jaeger}
Fran{\c{c}}ois Jaeger, \emph{On graphic-minimal spaces}, Ann. Discrete Math.
  \textbf{8} (1980), 123--126, Combinatorics 79 (Proc. Colloq., Univ.
  Montr\'eal, Montreal, Que., 1979), Part I.

\bibitem{JT}
Tommy~R. Jensen and Bjarne Toft, \emph{Graph coloring problems},
  Wiley-Interscience Series in Discrete Mathematics and Optimization, John
  Wiley \& Sons Inc., New York, 1995, A Wiley-Interscience Publication.

\bibitem{Kelmans}
Alexander~K. Kelmans, \emph{On edge bijections of graphs}, Tech. Report 93-41,
  DIMACS, 1993.

\bibitem{KR}
V{{\'a}}clav Koubek and Vojt{{\v e}}ch R{{\"o}}dl, \emph{On the minimum order
  of graphs with given semigroup}, J. Combin. Theory Ser. B \textbf{36} (1984),
  no.~2, 135--155.

\bibitem{Lang}
Serge Lang, \emph{Algebra}, Graduate Texts in Mathematics, vol. 211,
  Springer-Verlag, New York, 2002.

\bibitem{LMT}
Nathan Linial, Roy Meshulam, and Michael Tarsi, \emph{Matroidal bijections
  between graphs}, J. Combin. Theory Ser. B \textbf{45} (1988), no.~1, 31--44.

\bibitem{Reza}
Reza Naserasr, \emph{Homomorphisms and edge colourings of planar graphs}, (to
  appear in J. Combin. Theory Ser. B).

\bibitem{Nes-Ramsey}
Jaroslav Ne{\v{s}}et{\v{r}}il, \emph{Ramsey theory}, Handbook of combinatorics
  (R.L. Graham, M.~Gr{\"o}tschel, and L.~Lov{\'a}sz, eds.), Elsevier,
  Amsterdam, 1995, pp.~1331--1403.

\bibitem{NR-Ramsey}
Jaroslav Ne{\v{s}}et{\v{r}}il and Vojt{\v{e}}ch R{\"o}dl, \emph{Simple proof of
  the existence of restricted {R}amsey graphs by means of a partite
  construction}, Combinatorica \textbf{1} (1981), no.~2, 199--202.

\bibitem{NR-rigid}
Jaroslav Ne{\v{s}}et{\v{r}}il and Vojt{\v{e}}ch R{\"o}dl, \emph{Chromatically
  optimal rigid graphs}, J. Combin. Theory Ser. B \textbf{46} (1989), no.~2,
  133--141.

\bibitem{NT}
Jaroslav Ne{\v{s}}et{\v{r}}il and Claude Tardif, \emph{Duality theorems for
  finite structures (characterising gaps and good characterisations)}, J.
  Combin. Theory Ser. B \textbf{80} (2000), no.~1, 80--97.

\bibitem{NS}
Jaroslav Ne{\v s}et{\v r}il and Robert {{\v S}}{\'a}mal, \emph{On tension
  continuous mapings}, in preparation.

\bibitem{NR}
Jaroslav Ne{{\v s}}et{{\v r}}il and Vojt{{\v e}}ch R{\"o}dl, \emph{On {R}amsey
  graphs without cycles of short odd lengths}, Comment. Math. Univ. Carolin.
  \textbf{20} (1979), no.~3, 565--582.

\bibitem{RSthesis}
Robert {\v S}\'amal, \emph{On {XY} mappings}, Ph.D. thesis, Charles University,
  2005.

\bibitem{Wel}
Emo Welzl, \emph{Color-families are dense}, Theoret. Comput. Sci. \textbf{17}
  (1982), no.~1, 29--41.

\bibitem{Whit1}
Hassler Whitney, \emph{{Congruent graphs and the connectivity of graphs}}, Am.
  J. Math. \textbf{54} (1932), 150--168.

\bibitem{Whit2}
Hassler Whitney, \emph{{2-isomorphic graphs.}}, Am. J. Math. \textbf{55}
  (1933), 245--254.

\end{thebibliography}

\providecommand{\bysame}{\leavevmode\hbox to3em{\hrulefill}\thinspace}
\providecommand{\MR}{\relax\ifhmode\unskip\space\fi MR }
\providecommand{\MRhref}[2]{%
  \href{http://www.ams.org/mathscinet-getitem?mr=#1}{#2}
}
\providecommand{\href}[2]{#2}

\end{document}